\documentclass[11pt]{article}
\usepackage{amsmath}
\usepackage{graphicx}
\usepackage{enumerate}
\usepackage{url} 
\usepackage{epstopdf}
\usepackage{epsfig}
\usepackage{multirow}
\usepackage{graphics}
\usepackage{subfigure}
\usepackage{caption}
\usepackage{float}

\newtheorem{thm}{Theorem}
\newtheorem{defn}{Definition}

\begin{document}

\def\spacingset#1{\renewcommand{\baselinestretch}%
{#1}\small\normalsize} \spacingset{1}


  \title{\bf On the Assessment of Bootstrap Intervals for Samples of Fixed Size}
  \author{Weizhen Wang\hspace{.2cm}\\
  		School of Mathematics, Statistics and Mechanics, Beijing University of Technology\\
  		Department of Mathematics and Statistics, Wright State University\\ 
  		Chongxiu Yu \\
  		School of Mathematics, Statistics and Mechanics, Beijing University of Technology\\
  		and \\
  		Zhongzhan Zhang \\
  		School of Mathematics, Statistics and Mechanics, Beijing University of Technology
  		}

\date{2019/12/15}
\maketitle

\bigskip
\begin{abstract}
	 A reasonable confidence interval should have a confidence coefficient no less than the given nominal level $1-\alpha$ and a small expected length to reliably and accurately estimate the parameter of interest, and the bootstrap interval is considered to be an efficient interval estimation technique.
	 In this paper, we offer a first attempt at computing the coverage probability and expected length of a parametric or percentile bootstrap interval by exact probabilistic calculation for any fixed sample size. This method is applied to the basic bootstrap intervals for functions of binomial proportions and a normal mean. None of these intervals, however, are found to have a correct confidence coefficient, which leads to illogical conclusions including that the bootstrap interval is narrower than the $z$-interval when estimating a normal mean. This raises a general question of how to utilize bootstrap intervals appropriately in practice since the sample size is typically fixed. 

\end{abstract}

\noindent%

\vfill
\newpage
\spacingset{1.5} 
\section{Introduction}
Bootstrap confidence intervals, proposed by Efron (1979), are widely used in statistical practice. It uses a resampling technique to utilize the information repeatedly from the observed random sample and generate a randomized interval to estimate the parameter of interest. Efron, Rogosa and Tibshirani (2004) pointed out that even setting the number of bootstrap samples, the quantity $m$ to be introduced in (\ref{iid}), at 50 is likely to lead to fairly good standard error estimates.

There are three factors in evaluating a confidence interval: reliability, accuracy and accessibility.
The bootstrap interval is easy to access by a straightforward calculation.
How to assess the other two factors that are measured by coverage probability and expected length, respectively, for a sample of fixed size $n$? This is the goal of the paper. The current practice is mainly based on asymptotic analyses or statistical simulations. However, the obvious drawbacks
include no discussion on the case of not-large $n$ and limited simulations over a parameter space of infinite points. In fact, a negative signal on reliability was found in 
 Wang (2013). He proved that any bootstrap interval for a function of proportions, including the proportion, the difference and odds ratio, always has an infimum of coverage probability (ICP) zero (i.e., a zero confidence coefficient) regardless of $n$, the nominal level $1-\alpha$ and the estimator. In order not to cause ambiguity, the confidence coefficient of a confidence interval is defined to be the infimum of its coverage probability function over the entire parameter or distribution space (see Casella and Berger 2002, p. 418). 

This paper moves one step further to compute the coverage probability and expected length explicitly at any point in the parameter or distribution space for any given $n$. Then the reliability and accuracy of bootstrap interval can be checked directly. 
This approach is applied to
estimating a proportion, two functions of two proportions and a normal mean with or without a known variance $\sigma^2$.  We focus on two important cases: percentile bootstrap (Efron and Tibshirani, 1993, p. 171) and parametric bootstrap.  An effective probabilistic tool is provided to evaluate a bootstrap interval for any fixed $n$.  This is different from existing approaches (see, for examples, DiCiccio and Romano, 1988; DiCiccio and Tibshirani, 1987; Efron and Tibshirani, 1993; Shao and Tu, 1995; Mantalos and Zografos, 2008) that evaluate the bootstrap interval asymptotically or by simulations.  

In Section 2, we describe a general approach to calculate the coverage probability and expected length of parametric and percentile bootstrap intervals. This approach is applied to the bootstrap intervals for a proportion and two functions of two proportions in Sections 3 and 4, respectively. We then discuss bootstrap intervals for a normal mean with a known or unknown variance $\sigma^2$ in Sections 5 and 7, respectively. Section 6 deals with the estimation of a median. Section 8 contains discussions. All proofs are given in the Appendix.

\section{A general approach to compute the coverage probability and expected length of a bootstrap interval}

Suppose ${\cal F}$ is the distribution space which consists of a collection of cumulative distribution functions (CDF) $F(y)$. If all $F(y)$'s are determined by a parameter vector $\underline{\theta}=(\theta,\underline{\eta})$, then ${\cal F}$ is also called the parameter space and we have a parametric model. Here $\theta$ is the parameter of interest and $\underline{\eta}$ is the nuisance parameter vector.  If  $F(y)$'s are not determined by a parameter vector, 
we then have a nonparametric model. The goal is to estimate  $\theta$ 
under the parametric model or $\theta(F)$, a function of $F$, under the nonparametric model    
using $1-\alpha$ bootstrap intervals based on an independently identically distributed (i.i.d.) sample $\underline{y}=\{y_1,...,y_n\}$ from $F$. 

For reliability, 
we wish to find out whether the confidence coefficient of the bootstrap interval is no smaller
than or close to the given $1-\alpha$ for a fixed $n$ by providing a formula for the coverage probability.  For accuracy, a formula for the expected length 
is of interest.

\subsection{ The coverage probability and expected length of a $1-\alpha$ parametric bootstrap interval 
}

Suppose 
an i.i.d. sample  $\underline{y}=\{y_1,...,y_n\}$ of size $n$ is observed from 
a CDF $F_{\underline{\theta}}$ for $\underline{\theta}=(\theta,\underline{\eta})$. 
Let $\hat{\underline{\theta}}(\underline{y})=(\hat{\theta}(\underline{y}), \hat{\underline{\eta}}(\underline{y}))$ be the maximum likelihood estimator (MLE) for $\underline{\theta}$ over the parameter space $\cal F$ and let  $G_{\underline{\theta}}$ be the joint CDF of  $\hat{\underline{\theta}}$.

For  a fixed positive integer $m$ and  $\hat{\underline{\theta}}(\underline{y}) $ we observe a sample $\{\underline{y}_j^B=(y_{1j}^B,...,y_{nj}^B)\}_{j=1}^m$
of size $m$, where $\underline{y}_j^B$ for $j=1,...,m$ are i.i.d. and, for a fixed $j$,  $y_{ij}^B$ for $i=1,...,n$ are also i.i.d. and follow the CDF $F_{\hat{\underline{\theta}}(\underline{y})}$.  
We compute $\hat{\theta}(\underline{y}_j^B)$ for $j=1,...,m$. For the  observed $\hat{\underline{\theta}}(\underline{y})$, 
\begin{equation}
\label{iid}  \mbox{$\hat{\theta}(\underline{y}_j^B)|\hat{\underline{\theta}}(\underline{y})$ are i.i.d. for $j=1,...,m$ and follow a CDF,  $H_{\hat{\underline{\theta}}(\underline{y})}(x)$}.
\end{equation}
The CDF  $H_{\hat{\underline{\theta}}(\underline{y})}$ is determined by three items: the CDF $F_{\underline{\theta}}$, the estimator $\hat{\underline{\theta}}$ and the observation $\underline{y}$. In some commonly used parametric models, $H_{\hat{\underline{\theta}}(\underline{y})}$ can be derived explicitly as shown through  examples later.
The parametric bootstrap interval is given by two order statistics of the sequence $\{\hat{\theta}(\underline{y}_j^B)\}_{j=1}^m$. However, we need the following definition.

\begin{defn} For $p\in [0,1]$,
	the $100p$-th percentile of
	a data set $\{u_i\}_{i=1}^m$ of size $m$ is a quantity $k$ that at least $mp$ of $u_i$'s are no larger than $k$ and at least $m(1-p)$ of $u_i$'s are no smaller than $k$.
\end{defn}
Let $u_{(j)}$ denote the $j$th order statistic of $\{u_i\}_{i=1}^m$ for $j=1,...,m$. 
The $100p$-th percentile may or may not be unique depending on whether the product $mp$ is an integer. When $mp$ is not an integer, then $u_{([mp]+1)}$ is the unique $100p$-th percentile (here $[x]$ denotes the largest integer no larger than $x$); when $mp$ is an integer, then any value in interval $[u_{(mp)},u_{(mp+1)}]$ is a $100p$-th percentile, and $u_{(mp)}$ and  $u_{(mp+1)}$  are the smallest and largest $100p$-th percentiles, respectively. For examples, in a data set $\{u_i\}_{i=1}^m$ for $m=6$, any value in $[u_{(3)},u_{(4)}]$ is a $50$-th percentile for $p=0.5$; however, in a data set $\{u_i\}_{i=1}^5$, $u_{(3)}$ is the unique $50$-th percentile. In both cases, a percentile can be chosen to be an order statistic $u_{(j)}$ for some $j$, which is used for the entire paper. For simplicity,  we do not discuss the case that a percentile is equal to the average of two consecutive order statistics. 

Write $u_j=\hat{\theta}(\underline{y}_j^B)$. The  $1-\alpha$ parametric bootstrap interval for $\theta$ based on the set of bootstrap estimates $A_{\underline{y}}\stackrel{def}{=}\{u_j\}_{j=1}^m$ is 
\begin{equation}
\label{bootstrap-interval-p-g}
C_{pa}=
[u_{(m_l)},u_{(m_u)} ]\,\  \mbox{for $m_l=[m{\alpha\over 2}]+1$ and $m_u=m+1-m_l$},
\end{equation}
where $u_{(m_l)}$ is the largest  $100({\alpha\over 2})$-th percentile and  $u_{(m_u)}$ is the smallest  $100(1-{\alpha\over 2})$-th percentile of $A_{\underline{y}}$.
These choices make the interval (\ref{bootstrap-interval-p-g}) short.
e.g.,  
when $(m,\alpha)=(100, 0.1)$, $C_{pa}=
[u_{(6)},u_{(95)}]$. Both $u_{(5)}$ and $u_{(6)}$ are the 5-th percentile but we pick the larger one.

Let $Cover_{C_{pa}}(\theta,\underline{\eta})=P(u_{(m_l)}\leq \theta \leq u_{(m_u)} )$ be the coverage probability of $C_{pa}$. 
Then,
\begin{eqnarray*}
	\begin{split}
		Cover_{C_{pa}}(\theta,\underline{\eta})&\stackrel{*}{=}E_{\{\underline{y}\stackrel{iid}{\sim} F_{\underline{\theta}}\}} [P_{\underline{y}}(u_{(m_l)}\leq \theta \leq u_{(m_u)} |\underline{y})]
		 \stackrel{**}{=}E_{\{\hat{\underline{\theta}}(\underline{y})\sim G_{\underline{\theta} }\}} [P_{\hat{\underline{\theta}}(\underline{y})}(u_{(m_l)}\leq \theta \leq u_{(m_u)} |\hat{\underline{\theta}}(\underline{y}))]
	\end{split}
\end{eqnarray*}
From computational point of view, (**) is much simpler because it involves a lower-dimensional integration or summation,  a major advantage for the parametric model.  
Note
\begin{eqnarray*}
	P_{\hat{\underline{\theta}}}(u_{(m_l)}\leq \theta \leq u_{(m_u)} |\hat{\underline{\theta}})
	=	P_{\hat{\underline{\theta}} }(u_{(m_l)} \leq \theta |\hat{\underline{\theta}})-P_{\hat{\underline{\theta}} }(u_{(m_u)}<\theta |\hat{\underline{\theta}} ),
\end{eqnarray*}
and, as shown in Casella and Berger (2002, p. 229), the $j$-th order statistic $u_{(j)}$ has a CDF 
\begin{equation}
\label{order-CDF-p-g} 
F_j(x)
=\sum_{k=j}^m {m\choose k} H_{\hat{\underline{\theta}}}(x)^k(1-H_{\hat{\underline{\theta}}}(x))^{m-k}
=1-F_B(j-1,m,H_{\hat{\underline{\theta}}}(x))\,\ \mbox{for any}\,\  x\in R^1,
\end{equation} where $F_B(x,m,p)$ is the CDF of a binomial random variable $X\sim Bino(m,p)$. 
Thus,
\begin{equation}\label{cover-p-g-2}
Cover_{C_{pa}}(\theta,\underline{\eta})
=\int [F_B(m_u-1,m,H_{\hat{\underline{\theta}}}(\theta^-))-F_B(m_l-1,m,H_{\hat{\underline{\theta}}}(\theta))] G_{\underline{\theta}}(d\hat{\underline{\theta}}), 
\end{equation}
where $u_0^-$ denotes the (left) limit of $u$ when $u$ increasingly approaches $u_0$. We will use (\ref{cover-p-g-2}) repeatedly later. 
Following (\ref{order-CDF-p-g}), the expected length of $C_{pa}$, $E[u_{(mu)}-u_{(ml)}]$, is given by
\begin{equation}\label{el-p-g}
EL_{C_{pa}}(\theta,\underline{\eta})
=\int  [\int y F_B(m_u-1,m,H_{\hat{\underline{\theta}}}(dy))- \int y F_B(m_l-1,m,H_{\hat{\underline{\theta}}}(dy))] G_{\underline{\theta}}(d\hat{\underline{\theta}}).
\end{equation}

\subsection{The coverage probability and expected length of a $1-\alpha$ percentile bootstrap interval}
 
Under the nonparametric model, the goal is to estimate $\theta(F)$, a known function of the CDF $F$. We make it clear that a percentile interval can be derived in both nonparametric and parametric models.
Suppose an i.i.d. sample  $\underline{y}=\{y_1,...,y_n\}$ is observed from  
$F$ and $\hat{\theta}(\underline{y})$ is an estimator for $\theta(F)$. 
Let
$\underline{y}^B=\{y_1^B,...,y_n^B\}$ be a bootstrap sample with replacement 
from $\underline{y}$. If all $y_i$'s 
are distinct, then the induced sample space and probability mass are
\begin{equation}
\label{boot-s}
S_{\underline{y}}=\{\underline{y}^B=(y_1^B,...,y_n^B): y_j^B\in \underline{y}\}
\,\ 
\mbox{and}\,\
p_{\underline{y}}(\underline{y}^B)={1\over n^n},
\end{equation} respectively.
We observe an independent sample $\underline{y}_j^B=\{y_{1j}^B,...,y_{nj}^B\}_{j=1}^m$  of size $m$ from $S_{\underline{y}}$. Then, for the given $\underline{y}$, $\hat{\theta}(\underline{y}_j^B)\mid \underline{y}$ are i.i.d. 
for $j=1,...,m$, and for each $j$
\begin{equation}
\label{fy}   
\hat{\theta}(\underline{y}_j^B)\mid \underline{y}\sim H_{\underline{y}}(x)\,\ \mbox{for a CDF}\,\
H_{\underline{y}}(x)={\#\,\ \{\underline{y}^B\in S_{\underline{y}}: \hat{\theta}(\underline{y}^B)\leq x\}\over n^n}\,\  \mbox{for any}\,\  x\in R^1,
\end{equation}
where $\#(A)$ is the number of elements of set $A$.
Different from $H_{\hat{\underline{\theta}}(\underline{y})}$ in (\ref{iid}),  
$H_{\underline{y}}$ is determined by $\hat{\underline{\theta}}$ and $\underline{y}$ but not $F$. 
If at least two $y_i$'s are identical,
(\ref{fy}) is still valid if one labels those identical $y_i$'s with different markers. 
For example, rewrite $\underline{y}=(0,0,1,1,1)$ as $(0_1,0_2,1_1,1_2,1_3)$. By doing so,  $S_{\underline{y}}$ still has $n^n$ elements and both $(\ref{boot-s})$ and (\ref{fy}) hold.

The  $1-\alpha$ percentile interval for $\theta(F)$ based on $A_{\underline{y}}=\{u_j\}_{j=1}^m\stackrel{def}{=}\{\hat{\theta}(\underline{y}_j^B)\}_{j=1}^m$ is
$$
C_{pe}=[u_{(m_l)},u_{(m_u)}]\,\ \mbox{for}\,\ m_l=[m{\alpha \over 2}]+1, \,\ \mbox{and} \,\ m_u=m+1-m_l,
$$
which is the same as $C_{pa}$ in (\ref{bootstrap-interval-p-g}).
Its coverage probability at a given CDF $F$ is
\begin{equation}\label{cover-n-g}
Cover_{C_{pe}}(F) = \int...\int [F_B(m_u-1,m,H_{\underline{y}}(\theta(F)^-))-
F_B(m_l-1,m,H_{\underline{y}}(\theta(F)) ] F(dy_1)...F(dy_n),
\end{equation}
 which is much more complicated than (\ref{cover-p-g-2}). Similar to (\ref{el-p-g}), the expected length of $C_{pe}$ is
\begin{equation}\label{el-np-g}
EL_{C_{pe}}(F) =\int...\int [\int y F_B(m_u-1,m,H_{\underline{y}}(dy))- \int y F_B(m_l-1,m,H_{\underline{y}}(dy))] F(dy_1)...F(dy_n). 
\end{equation}
In this section, we derive four formulas (\ref{cover-p-g-2}), (\ref{el-p-g}), (\ref{cover-n-g}) and (\ref{el-np-g}) for coverage probability and expected length.

\section{Two bootstrap intervals for a proportion based on an i.i.d. Bernoulli sample}

 An i.i.d. sample $\underline{y}=\{y_1,...,y_n\}$ is observed from a Bernoulli population, $Bino(1,p)$. Interval estimation of $p$ is one of the basic statistical inference problems. There are a plenty of confidence intervals for $p$ based on the total $y=\sum_{i=1}^n y_i\sim Bino(n,p)$. Here we name only a few.  Asymptotic intervals: Wald interval,  Wilson interval (1927), and  Agresti-Coull interval (1998), and exact intervals: Clopper-Pearson interval (1934), Blyth-Still interval (1983), Casella refinement (1986) and Wang interval (2014). Now we use both parametric and percentile bootstrap  intervals to estimate $p$. The point estimators in the three asymptotic intervals have a form of $ay+b$ for two constants $a>0$ and $b$. Let
\begin{equation}
\label{three-asymptotic}
\left\{ \begin{tabular}{lll}
$\hat{p}={1\over n}y$& & for Wald interval;\\
$\breve{p}={1\over n+z^2_{{\alpha\over 2}} } y+ { z^2_{{\alpha\over 2}}/ 2 \over n+z^2_{{\alpha\over 2}}}$& $\stackrel{def}{=}a_1y+b_1 $   & for Wilson interval;\\
$\tilde{p}={1\over n+z^2_{{\alpha\over 2}} } y+  { z^2_{{\alpha\over 2}}/ 2 \over n+z^2_{{\alpha\over 2}}}$ &   & for Agresti-Coull interval,
\end{tabular}
\right. 
\end{equation}
where $z_{\alpha\over 2}$ is the $(1-{\alpha \over 2})$-th percentile of the standard normal distribution. Note $\breve{p}=\tilde{p}$. The bootstrap intervals based on $\breve{p}$ and $\tilde{p}$ are identical. So, intervals are derived using $\hat{p}$ and $\breve{p}$.

\subsection{The parametric and percentile bootstrap intervals based on $\hat{p}$}
 Here are the steps to build the parametric bootstrap interval using $\hat{p}$. i) Generate an i.i.d. sample: $\{y_1,...,y_n\}\sim Bino(1,p)$. 
Compute $\hat{p}={y\over n } $, the MLE and the uniformly minimum variance unbiased estimator (UMVUE) of $p$. Then $y=n\hat{p}\sim Bino(n,p)$. ii) For this given $\hat{p}$, generate an i.i.d. sample:
$\{y_{11}^B,...,y_{n1}^B|\hat{p}\}\sim Bino(1,\hat{p}).$
Compute $\hat{p}^{B}_1={1\over n } \sum_{i=1}^{n}y_{i1}^B$. Then $n\hat{p}^{B}_1|\hat{p} \sim  Bino(n,{y\over n}).$ iii) Repeat step ii) for $m-1$ more times. For the given $\hat{p}$ obtain 
\begin{equation}
\label{distphat}
\{\hat{p}^{B}_1,...,\hat{p}^{B}_m|\hat{p}\}\stackrel{iid}{\sim}  {1\over n}Bino(n,{y \over n}).
\end{equation}
iv) The $1-\alpha$ parametric bootstrap interval for $p$ based on $\hat{p}$ is $C_{wa}=[\hat{p}_{(m_l)}^B, \hat{p}_{(m_u)}^B]$.

Using (\ref{cover-p-g-2}) here,  $G_{\underline{\theta}}$ and $H_{\hat{\underline{\theta}}}(\underline{y})$ are CDFs for $y\sim Bino(n,{y\over n})$ and $n\hat{p}^B_j
|y\sim Bino(n,{y\over n})$, respectively. Then the  coverage probability function for interval $C_{wa}$ is 
\begin{equation}\label{cp-wald}
Cover_{C_{wa}}(p)= \sum_{y=0}^n[ F_B(m_u-1,m,F_B((np)^-,n,{y\over n})) -F_B(m_l-1,m,F_B(np,n,{y\over n})] p_B(y,n,p),
\end{equation}
where $p_B(x,n,p)$ is the probability mass function (PMF) of $Bino(n,p)$.
Clearly, the ICP 
is zero and is achieved when $p$ goes to zero or one for any sample size $n$ and any level $1-\alpha$. This fact was
first found in Wang (2013) but is confirmed here by (\ref{cp-wald}). Furthermore, using (\ref{cp-wald}) we can obtain the coverage probability at any value of $p$ explicitly. 

We now consider the percentile bootstrap interval for $p$ based on $\hat{p}$. In the observed binary sample $\underline{y}=\{y_1,...,y_n\}$, recall $y=\sum_{i=1}^n y_i$ and $\hat{p}={y \over n }$. Then the bootstrap sample proportion $\hat{p}^B={\sum_{i=1}^n y_i^B\over n}$ satisfies
$n\hat{p}^B|\hat{p} \sim Bino(n,{y\over n}).$
Comparing this with Step ii) that generates $C_{wa}$, we find an identical distribution and then conclude that the percentile bootstrap interval in this case is equal to the parametric bootstrap interval $C_{wa}$. 

\subsection{The bootstrap interval for $p$ based on $\breve{p}$}

The construction of the $1-\alpha$ parametric bootstrap interval based on $\breve{p}$ is similar to $C_{wa}$. i.e., i) Generate an i.i.d. sample of size $n$: $\{y_1,...,y_n\}\sim Bino(1,p)$. 
Compute $\breve{p}$, the center of Wilson interval. Then $y={\breve{p}-b_1 \over a_1} \sim Bino(n,p)$. ii) For this given $\breve{p}$, generate an i.i.d. sample:
$\{y_{11}^B,...,y_{n1}^B|\breve{p}\} \sim Bino(1,\breve{p}).$
Compute $\breve{p}^{B}_1=a_1 \sum_{i=1}^{n}y_{i1}^B+b_1$. Then $ {\breve{p}^{B}_1-b_1\over a_1}|\breve{p} \sim  Bino(n,\breve{p}).$ iii) Repeat step ii) for $m-1$ more times and obtain 
\begin{equation}
\label{distwilson-s}
\{ {\breve{p}^{B}_1-b_1\over a_1},..., {\breve{p}^{B}_m-b_1\over a_1}|y\}\stackrel{iid}{\sim}  Bino(n,\breve{p})=Bino(n,a_1y+b_1).
\end{equation}
iv) The $1-\alpha$ parametric bootstrap interval based on  $\breve{p}$ is
$C_{wi}=[\breve{p}_{(m_l)}^B, \breve{p}_{(m_u)}^B].$

Following (\ref{cover-p-g-2}), the  coverage probability for $C_{wi}$ 
is 
\begin{equation}\label{cp-wilson}
\begin{split}
&Cover_{C_{wi}}(p)
= \sum_{y=0}^n[ F_B(m_u-1,m,F_B(({p-b_1\over a_1})^-,n,a_1y+b_1))\\& \hspace{1.3in} 
-F_B(m_l-1,m,F_B({p-b_1\over a_1},n,a_1y+b_1)] 
p_B(y,n,p).
\end{split}
\end{equation}
The ICP 
for  $C_{wi}$ is also zero  for any $n$ and $\alpha$. The $1-\alpha$ percentile interval based on $\breve{p}$ is identical to $C_{wi}$. In fact, the parametric and percentile bootstrap intervals based on the same estimator are identical
due to the binary sample $\underline{y}$ and the distribution $Bino(1,p)$.

\subsection{Interval comparisons}

We evaluate five intervals, two asymptotic intervals (Wald and Wilson), two bootstrap intervals $C_{wa}$ and $C_{wi}$, and Wang exact interval (2014), in terms of coverage probability and expected length. 
The expected length of $C_{wa}=[\hat{p}_{(m_l)}^B, \hat{p}_{(m_u)}^B]$, as shown in the Appendix, is
\begin{equation}
\label{el-wild} 
EL_{C_{wa}}(p)={1\over n}\sum_{y=0}^n \sum_{x=0}^{n-1}[F_B(m_u-1,m,F_B(x,n,{y\over n})) -F_B(m_l-1,m,F_B(x,n,{y\over n}))]p_B(y,n,p).
\end{equation}
The expected length of  $C_{wi}=[\breve{p}_{(m_l)}^B, \breve{p}_{(m_u)}^B]$ is derived similarly following (\ref{order-CDF-p-g}), (\ref{el-p-g}) and (\ref{distwilson-s}).
\begin{eqnarray*}
	EL_{C_{wi}}(p)=a_1\sum_{y=0}^n \sum_{x=0}^{n-1}[F_B(m_u-1,m,F_B(x,n,a_1y+b_1))-F_B(m_l-1,m,F_B(x,n,a_1y+b_1))]p_B(y,n,p).
\end{eqnarray*}

It is well known that Wald interval has a poor coverage probability since its confidence coefficient is equal to zero for any $n$ and $\alpha$. Huwang (1995) provided the asymptotic confidence coefficient for Wilson interval (1927).
e.g., when the nominal level $1-\alpha$ is 0.9, its  asymptotic confidence coefficient  is only 0.8000; 
while Wang exact interval  (2014) has a confidence coefficient 0.9 for any $n$.
\begin{figure}[h]
	\begin{center}
		\resizebox{0.8\textwidth}{0.29\textheight}
	{\begin{tabular}{ll}
		$C_{wa}, (1-\alpha,n,m)=(0.9,10,100)$ & $C_{wa}, (1-\alpha,n,m)=(0.9,30,10000)$ 	\\
		\epsfig{file=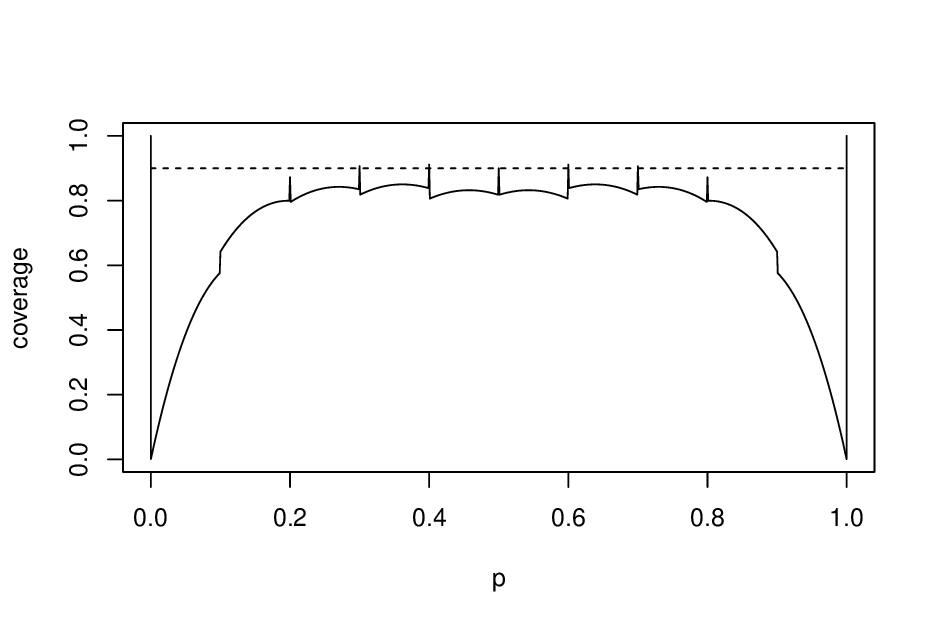,height=4cm, width=5.5cm,clip=}
		&  \epsfig{file=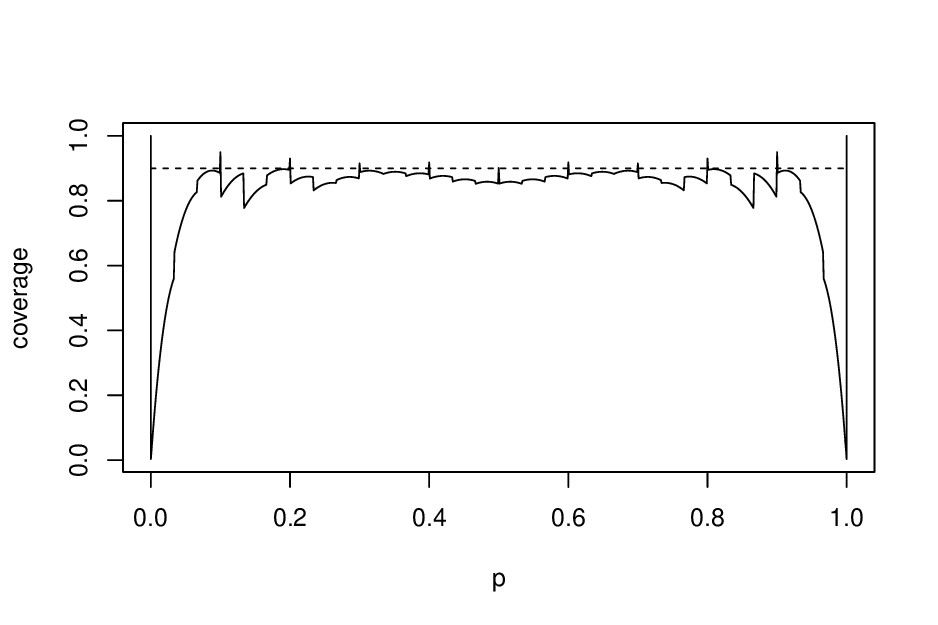,height=4cm,width=5.5cm,clip=} \\	
		$C_{wa}, (1-\alpha,n,m)=(0.9,200,100)$ &   $C_{wa}, (1-\alpha,n,m)=(0.9,200,10000)$	\\
		\epsfig{file=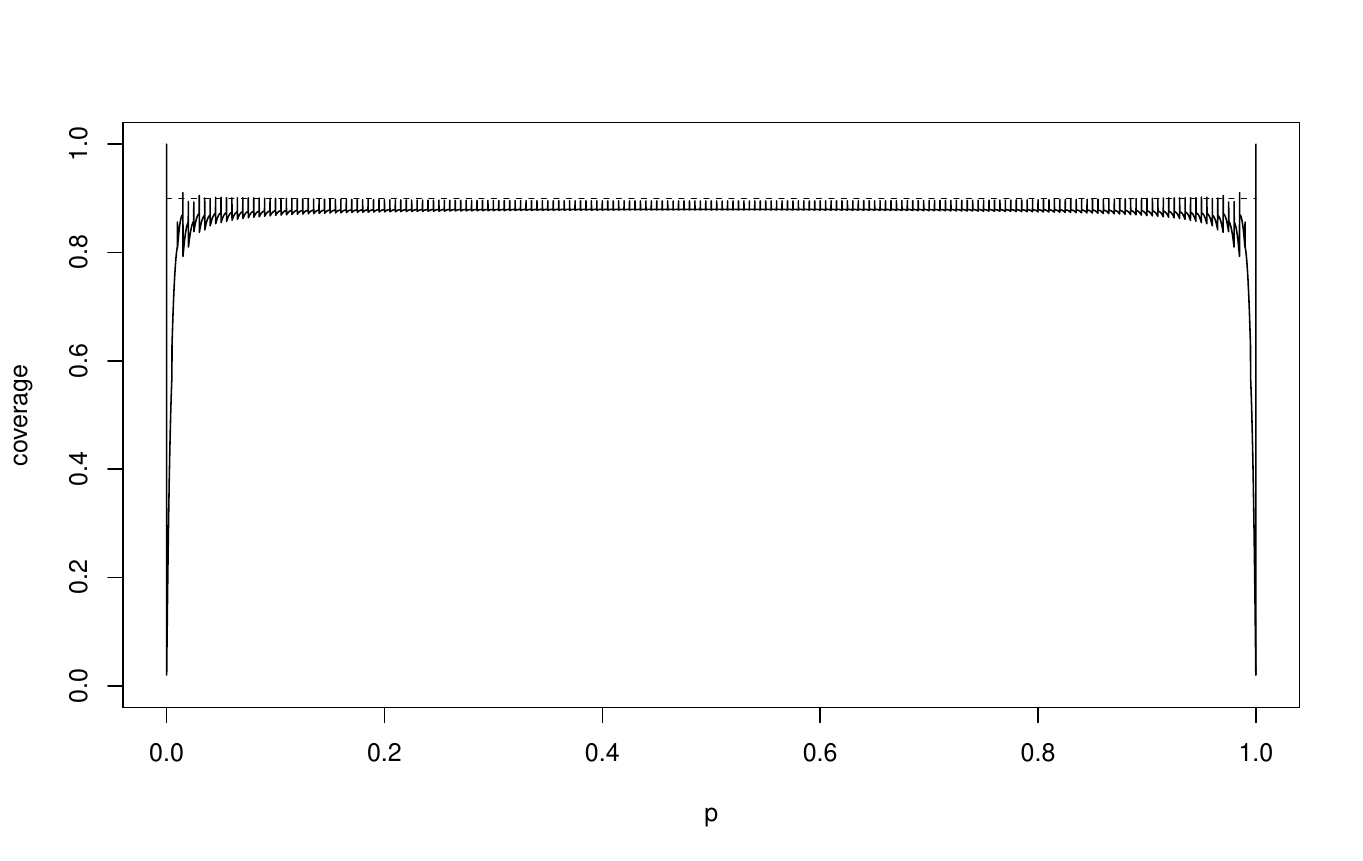,height=4cm, width=5.5cm,clip=}
		&  \epsfig{file=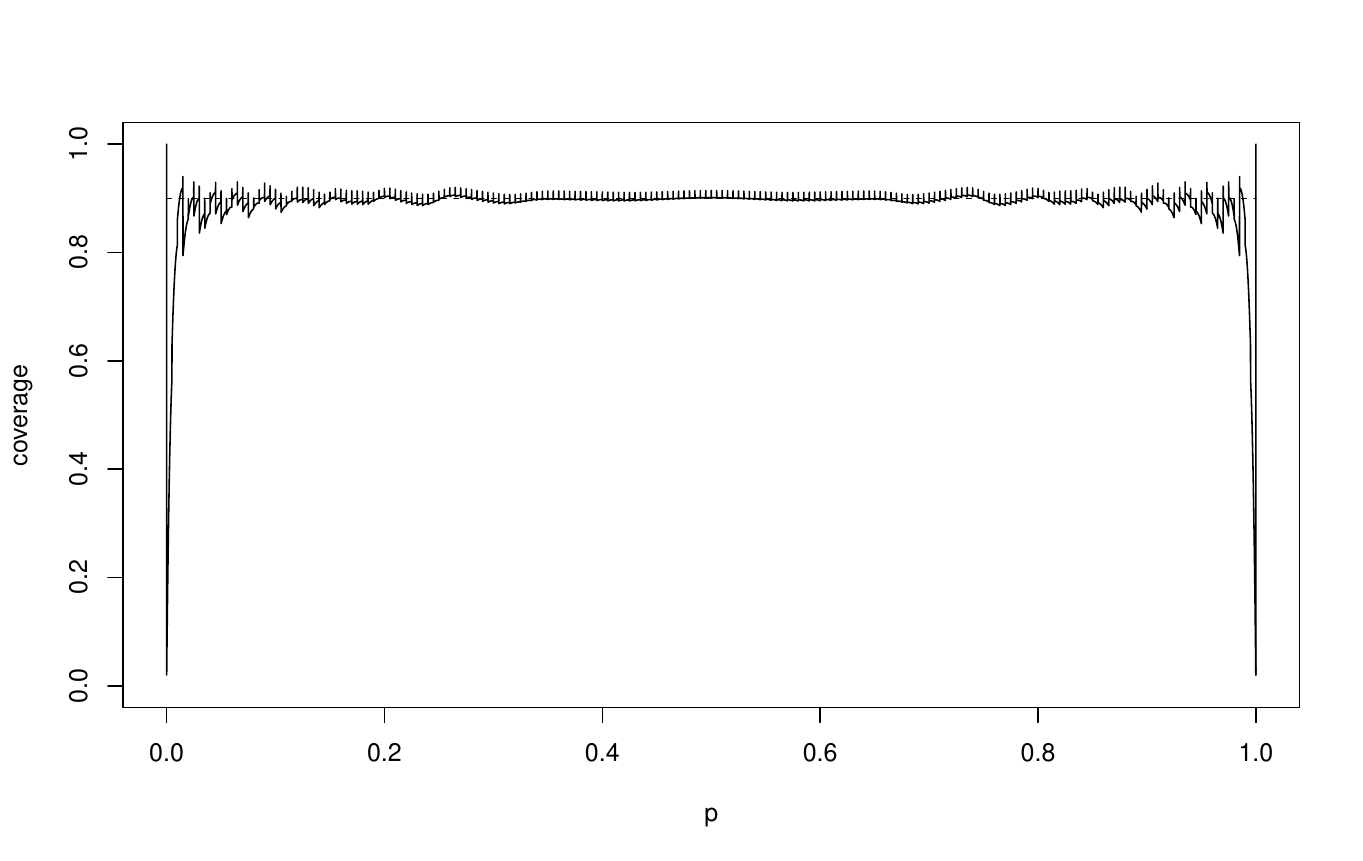,height=4cm,width=5.5cm,clip=} \\			
		$C_{wi}, (1-\alpha,n,m)=(0.9,10,100)$ & $C_{wi}, (1-\alpha,n,m)=(0.999,10,100)$ 	\\
		\epsfig{file=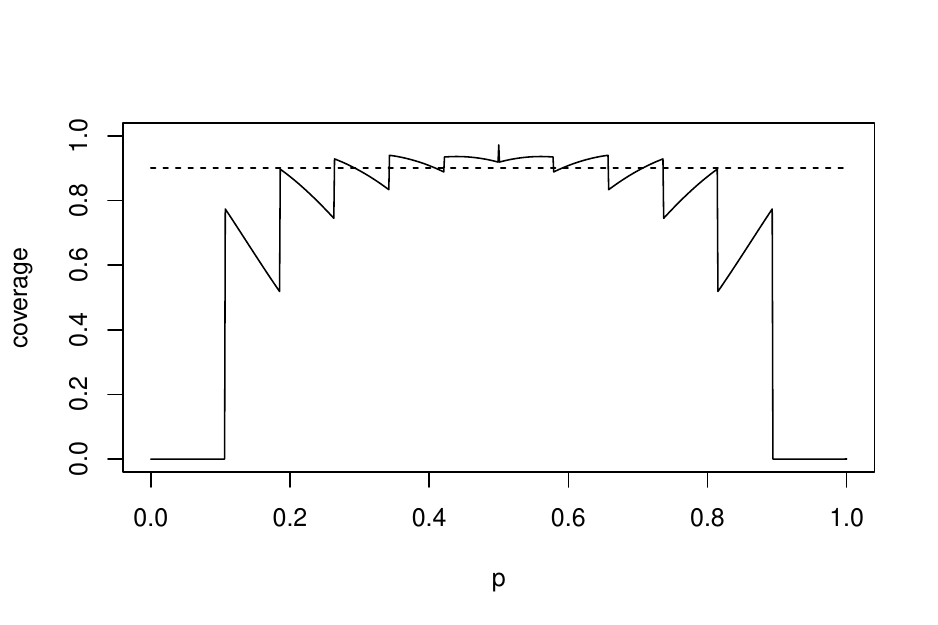,height=4cm, width=5.5cm,clip=}
		&   \epsfig{file=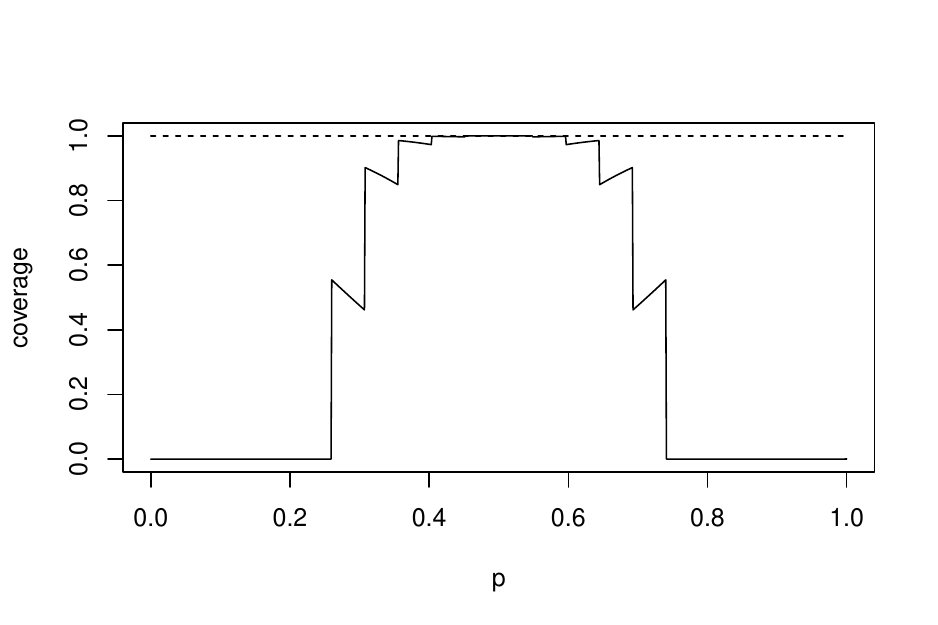,height=4cm,width=5.5cm,clip=} 		
	\end{tabular}}
	\caption{Coverage probabilities for $C_{wa}$ or $C_{wi}$ (solid) for different $(1-\alpha, n, m)$ and a reference line at $1-\alpha$ (dashed). }
	\label{fig1}
\end{center}
\end{figure}  
Figure~\ref{fig1} displays the coverage probabilities of the two bootstrap intervals $C_{wa}$ and $C_{wi}$. It coincides the findings of Wang (2013) that the confidence coefficient of any bootstrap interval for $p$ is always zero.  When $n$ or $m$ is not large, the coverage probability of the 90\% $C_{wa}$ is less than 0.9 for most values of $p$; when both $n$ and $m$ are large, it fluctuates around 0.9 for most values of $p$, however, the confidence coefficient is still zero. As seen in the second row of Figure~\ref{fig1}, $m$ has a big impact on the coverage probability of $C_{wa}$.
This makes the choice of $m$ difficult to have a correct coverage probability without the formula (\ref{cp-wald}). The choice of $m=50$ by Efron et al. (2004) seems too optimistic. Also, the third row of Figure~\ref{fig1} and the rows in Table \ref{tab-wilson} show
that the coverage probability for $C_{wi}$ does not increase in the nominal level $1-\alpha$ for fixed $(n, m)$, neither is the area under its coverage probability curve. This clearly creates chaos in choosing an appropriate nominal level.

\vspace{-0in}
\begin{table}[h]
	\caption{The area under the coverage probability curve of $C_{wi}$, $A$, when $1-\alpha$ changes.\label{tab-wilson}}
	\vspace{-0.3in}
	\begin{center}
		 \resizebox{0.6\textwidth}{0.13\textheight}
		 {\begin{tabular}{l|llllll}
			&\multicolumn{5}{|c}{$1-\alpha$} \\
			$(n,m)$	      & 0.6   &$0.8$  & $0.9$  & $0.95$ & $0.99$ &$0.999$\\ \hline
			$(10,100)$    & 0.5513& 0.6680 & 0.6625 & 0.6419 & 0.5555 & 0.4198   \\
			$(10,10000)$  & 0.5358& 0.6857 & 0.6682 & 0.6630 & 0.5737 & 0.4771   \\
			$(30,100)$    & 0.5749& 0.7273 & 0.7658 & 0.7755 & 0.7235 & 0.5733  \\
			$(30,10000)$  & 0.5733& 0.7504 & 0.7862 & 0.7830 & 0.7409 & 0.6781   \\
			$(100,100)$   & 0.5804& 0.7594 & 0.8307 & 0.8656 & 0.8552 & 0.7404   \\
			$(100,10000)$ & 0.5945& 0.7746 & 0.8475 & 0.8734 & 0.8758 & 0.8399
		\end{tabular}}	
	\end{center}
\end{table}

\begin{figure}
	\begin{center}
		\resizebox{0.8\textwidth}{0.29\textheight}{
		\begin{tabular}{cc} \vspace{-0.05in}
				$(n,m)=(10,100)$& $(n,m)=(10,10000)$\\ 
				\epsfig{file=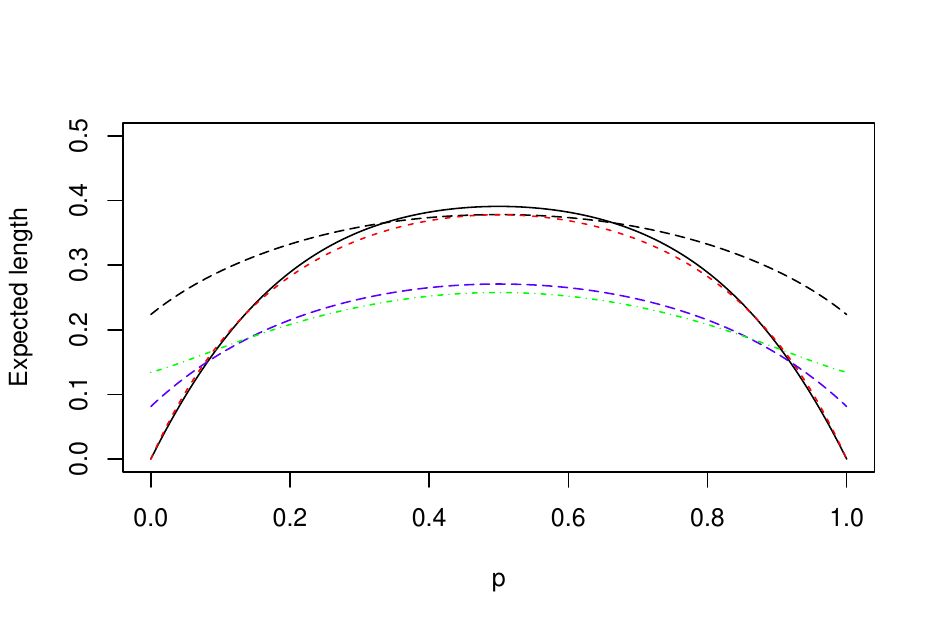,height=7.5cm, width=6cm,clip=}
				&  \epsfig{file=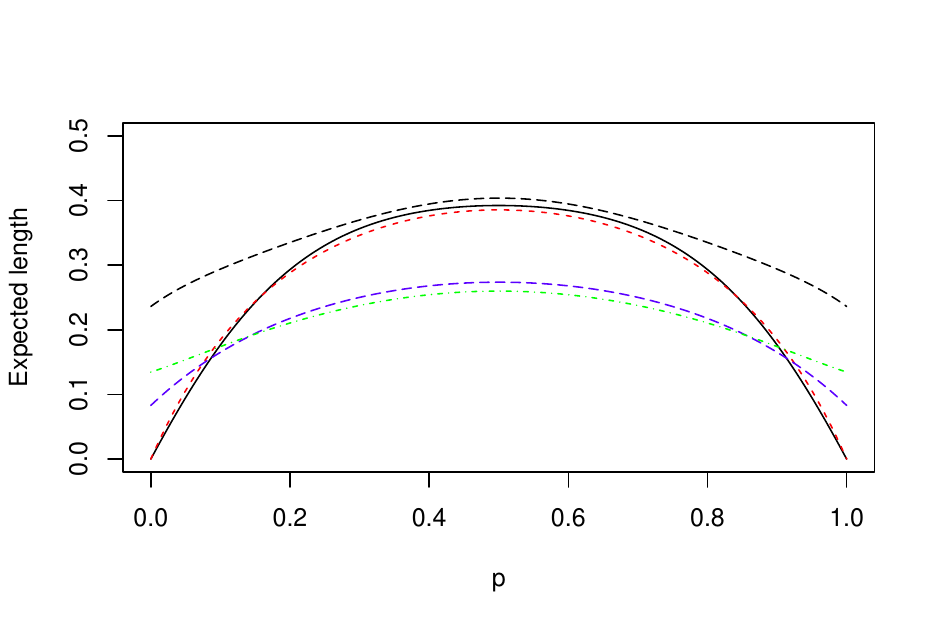,height=7.5cm,width=6cm,clip=} 
			 \\ \vspace{-0.05in}
				$(n,m)=(30,10000)$& $(n,m)=(100,100)$\\		
				\epsfig{file=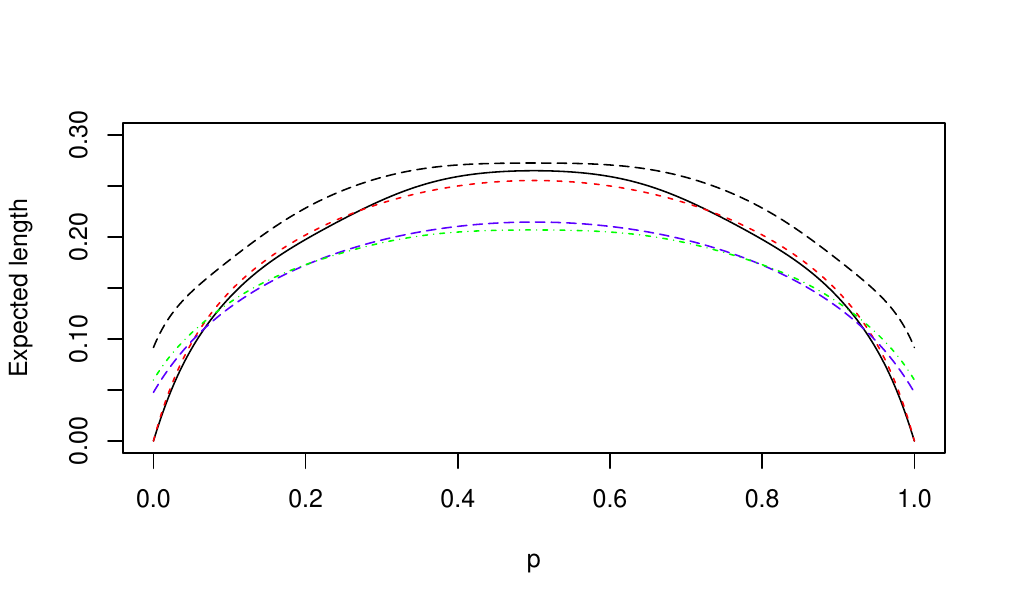,height=7.5cm, width=6cm,clip=}
				&  \epsfig{file=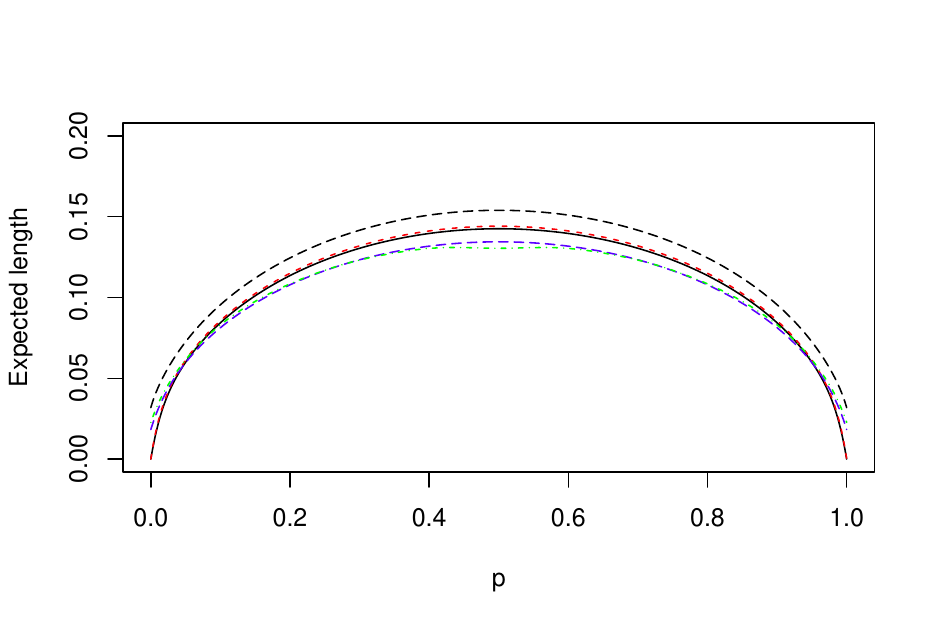,height=7.5cm,width=6cm,clip=}  \\
		\end{tabular}
	}
		\resizebox{0.9\textwidth}{0.16\textheight}{\begin{tabular}{l|l|lllll}
				\hline
				& 	      & $C_{wa}$   &Wald  & $C_{wi}$  & Wilson & Wang\\ \hline
				$(n,m)=$ & the nominal level $1-\alpha$  & 0.8050& 0.7942 & 0.9    & 0.6530 & 0.5265 \\
				(10,100)& confidence coefficient        & 0     & 0      & 0      & 0.3797 & 0.5265 \\
				A=0.6625& area under expected length& 0.2828& 0.2762 & 0.3328 & 0.2151 & 0.2132  \\ \hline
				$(n,m)=$ & the nominal level $1-\alpha$  & 0.7879& 0.8031 & 0.9    & 0.6585 & 0.5290 \\
				(10,10000)& confidence coefficient        & 0     & 0      & 0      & 0.3805 & 0.5289 \\
				A=0.6682& area under expected length& 0.2848& 0.2817 & 0.3440 & 0.2176 & 0.2153 \\ \hline
				$(n,m)=$ & the nominal level $1-\alpha$  & 0.8449& 0.8456 & 0.9    & 0.7797 & 0.7334 \\
				(30,10000)& confidence coefficient        & 0     & 0      & 0      & 0.5795 & 0.7334 \\
				A=0.7862& area under expected length& 0.1973& 0.1964 & 0.2262 & 0.1705 & 0.1703 \\ \hline
				$(n,m)=$ & the nominal level $1-\alpha$  & 0.8600& 0.8526 & 0.9    & 0.8274 & 0.8078 \\
				(100,100)& confidence coefficient        & 0     & 0      & 0      & 0.7263 & 0.8078 \\
				A=0.8307& area under expected length& 0.1111& 0.1125 & 0.1230 & 0.1061 & 0.1061 \\ 				
		\end{tabular}}
		\caption{Expected lengths for the five intervals with the same area under the coverage probability curve (A): $C_{wa}$ (black-solid), Wald (red), $C_{wi}$ (black-dash), Wilson (blue), and Wang (green), for different $n$ and $m$.\label{fig4}}
	\end{center}
\end{figure}

One common way of interval comparison is to compare the expected length of those intervals with the same confidence coefficient. This, however, cannot be done since the bootstrap intervals for $p$ always have a zero confidence coefficient. So, we compare the expected length of the intervals that have the same area under the coverage probability curve. 
Mantalos and Zografos (2008) did this comparison through limited simulations. Let 
$$A=\int_{0}^{1} Cover_{C_{wi}}(p)dp$$
for $Cover_{C_{wi}}(p).$ Then, we derive $C_{wa}$, Wald interval, Wilson interval, and Wang interval with the same area under the coverage probability curve, $A$, by choosing  a different
but appropriate value of $1-\alpha$ for each interval. This is doable because the coverage probability and expected length now can be computed precisely. The five expected lengths  are displayed in Figure~\ref{fig4}. In general, Wang interval has the smallest expected length but the highest confidence coefficient, while the two bootstrap intervals have largest expected lengths and a zero confidence coefficient. 
\section{Two bootstrap intervals for two functions of two proportions based on two i.i.d. Bernoulli samples}

It is often of interest to compare two treatments through the difference of two proportions, $d=p_1-p_2$, or the odds ratio, $\theta=p_1(1-p_2)/((1-p_1)p_2)$, using two independent Bernoulli samples $\{x_1,...,x_{n_1}\}\stackrel{iid}{\sim} Bino(1,p_1)$ and $\{y_1,...,y_{n_2}\}\stackrel{iid}{\sim} Bino(1,p_2)$. Let $\hat{p}_1=\sum_{i=1}^{n_1}x_i/n_1$ and $\hat{p}_2=\sum_{j=1}^{n_2}y_j/n_2$ be the two independent sample proportions. Then, $d$ and $\theta$ are estimated by $$\hat{d}=\hat{p}_1-\hat{p}_2\,\ \mbox{and} \,\  \hat{\theta}={(\hat{p}_1+{1\over 2n_1})(1-\hat{p}_2+{1\over 2n_2})\over (1-\hat{p}_1+{1\over 2n_1})(\hat{p}_2+{1\over 2n_2})},$$ respectively. For the estimation of $\theta$, we pick $\hat{\theta}$ by Gart (1966) rather than the commonly used MLE by Woolf (1955) because the MLE involves an undefined ratio $0/0$. We build a bootstrap interval for each of $d$ and $\theta$ based on $\hat{d}$ and $\hat{\theta}$. Different from Section 3, both cases involve a nuisance parameter, say $p_2$.
The interval construction is outlined below. 

i) Generate two i.i.d. samples: $\{x_1,...,x_{n_1}\}\sim Bino(1,p_1)$ and   $\{y_1,...,y_{n_2}\} \sim Bino(1,p_2)$.
Compute $\hat{p}_1$ and $\hat{p}_2$. Then $x=n_1\hat{p}_1\sim Bino(n_1,p_1)$ and $y=n_2\hat{p}_2\sim Bino(n_2,p_2)$. ii) For the given $\hat{p}_1$ and $\hat{p}_2$, generate two independent samples:
$\{x_{11}^B,...,x_{n_11}^B|\hat{p}_1\}\stackrel{iid}{\sim} Bino(1,\hat{p}_1),$ $ \{y_{11}^B,...,y_{n_21}^B|\hat{p}_2\}\stackrel{iid}{\sim} Bino(1,\hat{p}_2).$
Compute $\hat{p}^{B}_{11}= \sum_{i=1}^{n_1}x_{i1}^B/n_1$ and  $\hat{p}^{B}_{21}= \sum_{j=1}^{n_2}y_{j1}^B/n_2$. Then, $n_1\hat{p}^{B}_{11}|\hat{p}_1 \sim  Bino(n_1,{x\over n_1})$ and $n_2\hat{p}^{B}_{21}|\hat{p}_2 \sim  Bino(n_2,{y\over n_2})$ are independent. Compute
$$\hat{d}_1^B=\hat{p}^{B}_{11}-\hat{p}^{B}_{21}, \,\ \,\ \hat{\theta}_1^B={(\hat{p}^{B}_{11}+{1\over 2n_1})(1-\hat{p}^{B}_{21}+{1\over 2n_2})\over (1-\hat{p}^{B}_{11}+{1\over 2n_1})(\hat{p}^{B}_{21}+{1\over 2n_2})}. $$
iii) Repeat step ii) for $m-1$ more times and obtain 
\begin{equation}
\label{distdhat}
\{\hat{d}^{B}_1,...,\hat{d}^{B}_m|x,y\}\stackrel{iid}{\sim}  {1\over n_1}Bino(n_1,{x \over n_1})- {1\over n_2}Bino(n_2,{y \over n_2})
\end{equation}
with a CDF $H_{x,y}^D$ and 
$$\{\hat{\theta}^{B}_1,...,\hat{\theta}^{B}_m|x,y\}\stackrel{iid}{\sim}  {(Bino(n_1,{x \over n_1})+0.5)(n_2-Bino(n_2,{y \over n_2})+0.5)\over (n_1-Bino(n_1,{x \over n_1})+0.5)(Bino(n_2,{y \over n_2})+0.5)}
$$
with a CDF $H_{x,y}^R$. The two CDFs, $H_{x,y}^D$ and $H_{x,y}^R$, are given later.
iv) The two $1-\alpha$ parametric  bootstrap intervals for $d$ and $\theta$ 
are
$$C_{d}=[\hat{d}_{(m_l)}^B, \hat{d}_{(m_u)}^B]\,\  \mbox{and} \,\
C_{\theta}=[\hat{\theta}_{(m_l)}^B, \hat{\theta}_{(m_u)}^B],$$
respectively. The above intervals are also the percentile bootstrap intervals for $d$ and $\theta$. 

For the interval $C_{d}$, following (\ref{distdhat}), the conditional CDF of $\hat{d}_j^B|x,y$ is 
\begin{eqnarray*}
	H_{x,y}^D(z)
	&=&\sum_{v=0}^{n_2} F_B(n_1(z+{v\over n_2}),n_1,{x\over n_1})p_B(v,n_2,{y\over n_2}) \,\
	\mbox{for any} \,\   z\in R^1.
\end{eqnarray*} 
Then by (\ref{cover-p-g-2}) the coverage probability function for $C_d$ is shown to be
\begin{eqnarray*}
	\begin{split}
		&Cover_{C_d}(d,p_2)
		= \sum_{x=0}^{n_1}\sum_{y=0}^{n_2} [F_B(m_u-1,m,H_{x,y}^D(d^-))\\  &\hspace{1.15in}-F_B(m_l-1,m,H_{x,y}^D(d))]p_B(x,n_1,d+p_2)p_B(y,n_2,p_2).
	\end{split}
\end{eqnarray*}
This function is defined on the parameter space
$H_D=\{(d,p_2): p_2\in D(d)\,\ \mbox{for}\,\ d\in [-1,1] \},$
a parallelogram with  vertices (0,0), (1,0), (0,1) and (-1,1), and
$D(d)= [0,1-d]$ if $d\in [0,1]$; $D(d)=[-d,1]$ if $d\in [-1,0)$.
Following (\ref{order-CDF-p-g}), the conditional CDF of the $j$-th order statistic $\hat{d}_{(j)}^B\mid x,y$ is
$$F_j(z) =1-F_B(j-1,m,H_{x,y}^D(z))\,\ \mbox{for any} \,\ z\in R^1.$$
For given $x$ and $y$, $\hat{d}_j^B$ assumes values $\{d_s\}_{s=1}^K$ from the smallest $d_1$ to the largest $d_K$. Then, 
\begin{eqnarray}\label{ecd}
	E(\hat{d}^B_{(j)}|x,y) = d_1+ \sum_{s=1}^{K-1}(d_{s+1}-d_s)F_B(j-1,m,H_{x,y}^D(d_s)).
\end{eqnarray}
as shown in the Appendix. Therefore, the expected length of $C_{d}$ is given by
\begin{eqnarray*}
	\begin{split}
		EL_{C_{d}}(d,p_2)
		=\sum_{x=0}^{n_1}\sum_{y=0}^{n_2}&\{\sum_{s=1}^{K-1}(d_{s+1}-d_s)[F_B(m_u-1,m,H_{x,y}^D(d_s))\\ &\hspace{0.2in}-F_B(m_l-1,m,H_{x,y}^D(d_s))]\}p_B(x,n_1,d+p_2)p_B(y,n_2,p_2). 
	\end{split}
\end{eqnarray*}

For the interval $C_{\theta}$, the conditional CDF of $\hat{\theta}_j^B\mid x,y$ is 
\begin{eqnarray*}
		H_{x,y}^R(z)&=&
	\sum_{v=0}^{n_2} F_B( R(z,v),n_1,{x\over n_1})p_B(v,n_2,{y\over n_2})
	\,\ \mbox{for}\,\\
	 R(z,v)&=& {z(n_1+0.5)(v+0.5)-0.5(n_2-v+0.5) \over (n_2-v+0.5)+z(v+0.5) }.
\end{eqnarray*}
Then,  similar to the case of $C_d$, we have the following coverage probability function for $C_{\theta}$,
\begin{eqnarray*}
	\begin{split}
		Cover_{C_{\theta}}(\theta,p_2)
		&= \sum_{y=0}^{n_2}\sum_{x=0}^{n_1} [F_B(m_u-1,m,H_{x,y}^R(\theta^-))\\
		& \hspace{0.6in}-F_B(m_l-1,m,H_{x,y}^R(\theta))]p_B(x,n_1,{\theta p_2\over 1+(\theta-1)p_2})p_B(y,n_2,p_2).\\
	\end{split}
\end{eqnarray*}
on a domain of $H_R =\{(\theta,p_2): p_2\in [0,1],\theta\in[0,+\infty) \}$. 
The expected length for $C_{\theta}$ is
\begin{eqnarray*} 
	\begin{split}
		EL_{C_{\theta}}(\theta,p_2) 
		&=\sum_{y=0}^{n_2}\sum_{x=0}^{n_1}\{\sum_{s=1}^{K-1}(\theta_{s+1}-\theta_s)[F_B(m_u-1,m,H_{x,y}^R(\theta_s))\\ & \hspace{0.6in} -F_B(m_l-1,m,H_{x,y}^R(\theta_s))]\}p_B(x,n_1,{\theta p_2\over 1+(\theta-1)p_2})p_B(y,n_2,p_2), 
	\end{split}
\end{eqnarray*}
where, for given $x$ and $y$, 
 $\{\theta_s\}_{s=1}^K$ are the possible values of $\hat{\theta}_j^B$ arranged  from the smallest to the largest.

\begin{figure}	
	\subfigure[The 3-D coverage probability for $C_{d}$]{
		\begin{minipage}{0.5\textwidth}
			\centering
			\includegraphics[height=0.2\textheight,width=0.7\textwidth]{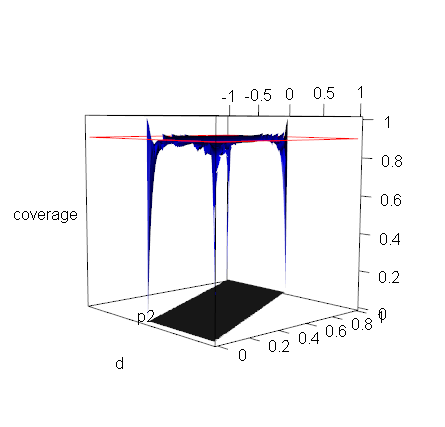}
			\includegraphics[height=0.2\textheight,width=0.7\textwidth]{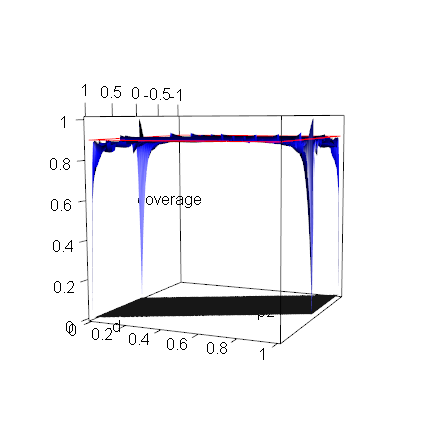}
			\includegraphics[height=0.2\textheight,width=0.7\textwidth]{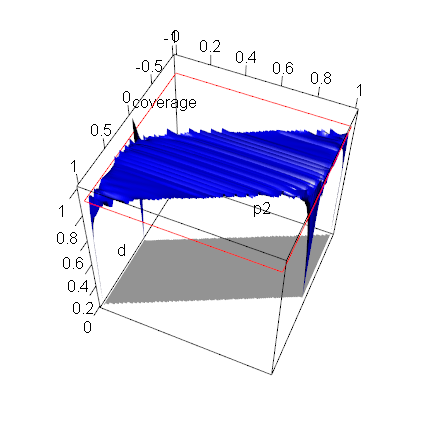}
		\end{minipage}
	}
	\subfigure[The 3-D coverage probability for $C_{\theta}$]{
		\begin{minipage}{0.5\textwidth}
			\centering
			\includegraphics[height=0.2\textheight,width=0.7\textwidth]{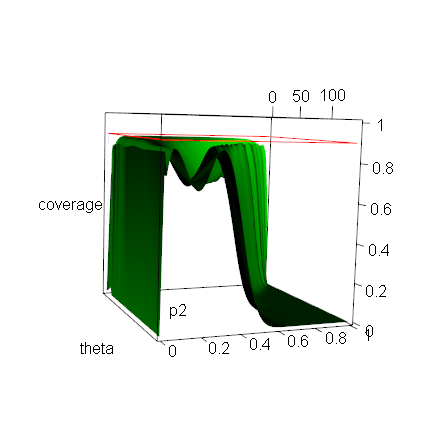}
			\includegraphics[height=0.2\textheight,width=0.7\textwidth]{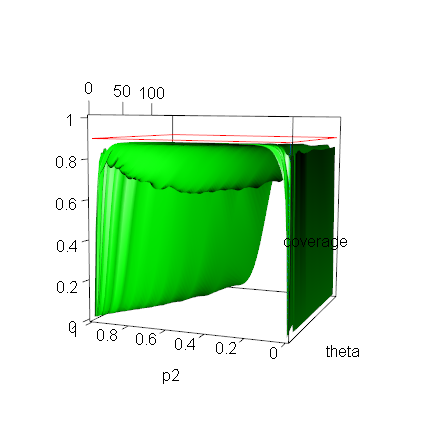}
			\includegraphics[height=0.2\textheight,width=0.7\textwidth]{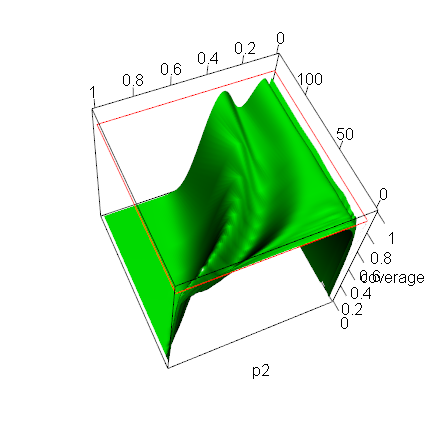}
		\end{minipage}
	}
	\caption{The 3-D plots of coverage probability for the 90\% bootstrap intervals: $C_d$ and $C_{\theta}$, when $(n_1,n_2,m)=(30, 60, 1000)$. The left column contains the coverage probability for $C_d$ (blue) on $H_{D}$ (the black parallelogram) and the plane $Coverage=0.9$ (the red square) at three different angles. The right column is the coverage probability for $C_{\theta}$ (green) as a function of $\theta \in [1,100]$ and $p_2\in [0,1]$ and the plane $Coverage=0.9$ (red) at three different angles. The ICPs
		for the two intervals are equal to 0. Most of the coverage 
		probability values are less than 0.9, especially for $C_{\theta}$.
	\label{fig3}} 
\end{figure}


Figure~\ref{fig3} displays the coverage probability functions of the 90\% bootstrap intervals, $C_d$ and $C_{\theta}$, for $(n_1,n_2,m)=(30,60,1000)$. If the intervals were truly of level 0.9, then the coverage probability surfaces would be on top of the plane $Coverage =0.9$, which, however, do not happen. In fact, their ICPs are both equal to zero, indicating that they are 0\% confidence intervals;  most of the coverage probability values are less than 0.9; between the two intervals, $C_{\theta}$ is more liberal.

\section{The parametric and percentile bootstrap intervals of a normal mean $\mu$ when the variance $\sigma^2$ is known}

An i.i.d. sample $\underline{y}=\{y_1,...,y_n\}$ is observed from a normal population $N(\mu,\sigma^2)$ with a probability density function (PDF) $p_N(x,\mu,\sigma)$ and a CDF $F_N(x,\mu,\sigma)$.  We wish to derive bootstrap intervals for $\mu$ based on an point estimator $Q=Q(\underline{y})$ that has a special PDF 
\begin{equation}\label{f-form}
p(q,\mu,\sigma)={1\over \sigma}f({q-\mu\over \sigma}) \,\ \mbox{for some even function
	$f(x)$}
\end{equation}
and a CDF $F(q,\mu,\sigma)$. Then, $(Q-\mu)/\sigma$ has a PDF
$p(x,0,1)=f(x)$ and a CDF $F(x,0,1)$. 

The first example of $Q$ is the sample mean, $\bar{y}$, since $\bar{y}$ has a PDF  in (\ref{f-form}) for $f(x)=p_N(x,0,{1\over \sqrt{n}})$.
The second example of $Q$ is the sample median. Let $a=[n/2]+1$. When $n$ is odd, then $Q=y_{(a)}$; otherwise,  $Q=(y_{(a-1)}+y_{(a)})/2$. For simplicity, we discuss an odd $n$ in this section. Then $y_{(a)}$ has a PDF in (\ref{f-form}) for $$f(x)=n\,\ p_N(x,0,1)p_B(a-1,n-1,F_N(x,0,1)).$$   
This function is even because $p_N(x,0,1)$ is even and 
$n=2a-1$ for an odd $n$. When $n$ is even, one can also show that the PDF of $Q$ satisfies (\ref{f-form}).  
Besides the two examples, $Q$ can be the truncated mean estimator for $\mu$ as well.

\subsection{The parametric bootstrap interval for $\mu$ based on $Q$}

The parametric bootstrap interval for $\mu$ based on $Q$ is constructed as follows.
i) Generate an i.i.d. sample, $\{y_1,...,y_n\}\sim N(\mu,\sigma^2)$. 
Compute $Q$. 
ii) For this $Q$, generate an i.i.d. bootstrap sample, $\underline{y}^B_1=\{y_{11}^B,...,y_{n1}^B\mid Q\} \sim N(Q,\sigma^2)$.
Compute $Q^{B}_1=Q(\underline{y}^B_1)$. Then $Q^{B}_1\mid Q \sim p(q,Q,\sigma).$
iii) Repeat ii) for $m-1$ more times and obtain 
$\{Q^{B}_1,...,Q^{B}_m\mid Q\}\stackrel{iid}{\sim} p(q,Q,\sigma).$
iv) The parametric bootstrap interval for $\mu$ based on $Q$ is
$C_Q=[Q_{(m_l)}^B, Q_{(m_u)}^B].$

\begin{thm}\label{thm-1}
	Let $Cover_{C_Q}(\mu)$ and  $EL_{C_Q}(\mu)$ be the coverage probability  and the expected length for $C_Q$, respectively. Then,
	\begin{equation}\label{cp-normal-q}
	Cover_{C_Q}(\mu) ={m_u-m_l\over m+1} 
	\end{equation}
	and
	\begin{equation}\label{el-normal-q}
	EL_{C_Q}(\mu) = m \sigma\int_{0}^{1} F^{-1}(z,0,1) [
	p_B(m_u-1,m-1,z) -p_B(m_l-1,m-1,z)]dz.
	\end{equation}
\end{thm}

Similar to the well-known $z$-interval,  $\bar{y}\pm z_{{\alpha\over 2}}\sigma/\sqrt{n}$,  the coverage probability $Cover_{C_Q}$ is independent of $(\mu,\sigma^2,n)$ but depends on $(m,\alpha)$. So, its confidence coefficient is equal to the constant coverage probability. 
$EL_{C_Q}(\mu)$ depends on $(\sigma,n,m,\alpha)$ but not $\mu$.	

\subsection{The parametric and percentile bootstrap intervals for $\mu$ based on the sample mean} 
When Theorem~\ref{thm-1} is applied to $Q=\bar{y}$, then the parametric bootstrap interval is
$C_{N}=[\bar{y}_{(m_l)}^B, \bar{y}_{(m_u)}^B]$
with the constant coverage probability in (\ref{cp-normal-q}) and the expected length	
\begin{equation}\label{elnormal}
\begin{split}
&EL_{C_N}(\mu)=\{{m\over \sqrt{n}}\int_{0}^{1} F_N^{-1}(z,0,1) [
p_B(m_u-1,m-1,z)\\
&\hspace{2.2in} -p_B(m_l-1,m-1,z)]dz \} \sigma \stackrel{def}{=}A(n,m,\alpha)\sigma
\end{split}
\end{equation}

\begin{table}
	\caption{The coverage probabilities ($CP$) and expected lengths ($EL$) of three sample-mean-based bootstrap intervals, $C_N$, $C_{pN}$ and $C_{Nu}$,   two sample-median-based bootstrap intervals, $C_{NM}$ and $C_{pM}$, and two classic confidence intervals,   
 the $z$-interval and the $t$-interval, $\bar{y}\pm t_{{\alpha\over 2},n-1}s/\sqrt{n}$,  for $\sigma=1$ and different $(n,m,1-\alpha)$. Each pair contains $(CP,EL)$.\label{tab-normal}}
	\begin{center}\vspace{-0.2in}
		\resizebox{0.8\textwidth}{0.4\textheight}{
			\begin{tabular}{ll|llll}
			\hline
         &	&$m=50$ & $m=100$ & $m=5000$ & $z$-interval$^+$ or $t$-interval$^-$\\
			\hline
		 &	&\multicolumn{4}{c}{$1-\alpha=0.9$}\\
			$n=5$ & $C_N$ & (0.8824, 1.4567) &  (0.8812, 1.4232) & (0.8996, 1.4702) 
			& (0.9, 1.4712)$^+$\\
			& $C_{pN}$ & (0.7707, 1.2232) & (0.7681, 1.1990) & (0.7860, 1.2453) & \\
			& $C_{Nu}$ & (0.7708, 1.2247) & (0.7671, 1.1966) &
			(0.7844, 1.2361) & (0.9, 1.7923)$^-$\\
			& $C_{NM}$ & (0.8824, 1.7444) & (0.8812, 1.7037) &
			(0.8996, 1.7601) & \\ 
			& $C_{pM}$ & (0.7998, 1.7372) & (0.7888, 1.6904) &
			(0.9350, 2.3153) & \\
			$n=31$& $C_N$   & (0.8824, 0.5850) &  (0.8812, 0.5716) & (0.8996, 0.5904) 
			& (0.9, 0.5908)$^+$\\
			& $C_{pN}$ & (0.8665, 0.5706) & (0.8654, 0.5581) & (0.8843, 0.5760)\\
			& $C_{Nu}$ & (0.8668, 0.5707) & (0.8651, 0.5576) &
			(0.8835, 0.5760) & (0.9, 0.6046)$^-$\\
			& $C_{NM}$   & (0.8824, 0.7281) &  (0.8812, 0.7113) & (0.8996, 0.7348) 
			 \\ 
			& $C_{pM}$   & (0.8721, 0.7229) &  (0.8705, 0.7060) & (0.8686, 0.6935) & \\					
			$n=301$ & $C_N$  & (0.8824, 0.1877) &  (0.8812, 0.1834) & (0.8996, 0.1895) & (0.9, 0.1896)$^+$\\
			& $C_{pN}$ & (0.8811, 0.1872) & (0.8804, 0.1831) & (0.8983, 0.1890)\\
			& $C_{Nu}$ & (0.8808, 0.1873) & (0.8796, 0.1830) & (0.8980, 0.1890) & (0.9, 0.1900)$^-$\\ 
			& $C_{NM}$  & (0.8824, 0.2351) &  (0.8812, 0.2297) & (0.8996, 0.2373) &\\
			& $C_{pM}$   & (0.8813, 0.2350) &  (0.8801, 0.2296) & (0.8966, 0.2357) & \\						
			&&\multicolumn{4}{c}{$1-\alpha=0.95$}\\
			$n=5$& $C_N$ & (0.9216, 1.6590) &  (0.9406, 1.7409) & (0.9496, 1.7513) 
			& (0.95, 1.7530)$^+$\\
			& $C_{pN}$ & (0.8096, 1.3729) & (0.8311, 1.4348) & (0.8378, 1.4490) \\
			& $C_{Nu}$ & (0.8163, 1.3948) & (0.8372, 1.4636) &(0.8451, 1.4724) 
			& (0.95, 2.3343)$^-$\\ 
			& $C_{NM}$    & (0.9216, 1.9891) &  (0.9406, 2.0878) & (0.9496, 2.0999) & \\
			& $C_{pM}$   & (0.8730, 2.0504) &  (0.9168, 2.2374) & (0.9375,    2.3259) & \\
			$n=31$ & $C_N$   & (0.9216, 0.6663) &  (0.9406, 0.6992) & (0.9496, 0.7033) 
			& (0.95, 0.7040)$^+$\\
			& $C_{pN}$ & (0.9084, 0.6489) & (0.9284, 0.6816) & (0.9356, 0.6861) &\\
			& $C_{Nu}$ & (0.9077, 0.6500) & (0.9273, 0.6821) & (0.9362, 0.6862) & (0.95, 0.7275)$^-$\\ 	
			& $C_{NM}$   & (0.9216, 0.8295) &  (0.9406, 0.8704) & (0.9496, 0.8756) & \\
			& $C_{pM}$   & (0.9124, 0.8235) &  (0.9318, 0.8642) & (0.9294, 0.8265) &\\									
			$n=301$ & $C_N$  & (0.9216, 0.2138) &  (0.9406, 0.2244) & (0.9496, 0.2257) & (0.95, 0.2259)$^+$\\
			& $C_{pN}$ & (0.9205, 0.2135) & (0.9395,0.2238) & (0.9462,0.2252)\\
			& $C_{Nu}$ & (0.9202, 0.2133) & (0.9393, 0.2238) & (0.9483, 0.2252) & (0.95, 0.2267)$^-$\\ 
			& $C_{NM}$ & (0.9216, 0.2678) & (0.9406, 0.2811) & (0.9496, 0.2827) &\\
			& $C_{pM}$ & (0.9207, 0.2676) & (0.9397, 0.2808) & (0.9491, 0.2829) &\\ 		  
		\end{tabular}
	}	
	\end{center}
\end{table}	
\noindent because $F^{-1}(z,0,1)=F_N^{-1}(z,0,1)/\sqrt{n}$.	

Table~\ref{tab-normal}  
contains exact calculation results using (\ref{cp-normal-q}) and (\ref{elnormal}). Three concerns are worth mentioning: a)  The coverage probabilities for $C_N$ are smaller than the nominal level $1-\alpha$ but approach $1-\alpha$ as $m$ goes large, much larger than 50. b) The coverage probability may not necessarily increase in $m$, e.g., the cases of $m=50$ and 100 when $1-\alpha=0.9$.  So, a larger $m$, i.e., more efforts, may not yield a more reliable result. 
c) Between the bootstrap interval $C_N$ and the $z$-interval, the former always has a shorter expected length than the latter. If deriving confidence intervals by a predetermined level $1-\alpha$, a typical approach in practice, one would conclude that the bootstrap interval $C_N$ dominates the $z$-interval in length even under the normal model -- an illogical conclusion that misleads practitioners. 


The percentile bootstrap interval for $\mu$ based on $\bar{y}$ is built as follows. 
i) Generate an i.i.d. sample: $\underline{y}=\{y_1,...,y_n\}\sim N(\mu,\sigma^2)$. 
Without loss of generality, we assume all $y_i$'s are distinct because this occurs with probability one.
ii) For this $\underline{y}$, consider $S_{\underline{y}}$ in (\ref{boot-s}). Define a random variable $\bar{y}^B=(y_1^B+...+y_n^B)/n$ on $S_{\underline{y}}$ with a CDF $H_{\underline{y}}$. Generate an i.i.d. bootstrap sample $\{y_{11}^B,...,y_{n1}^B\}$ from $\underline{y}$ 
and compute $\bar{y}_1^B=(y_{11}^B+...+y_{n1}^B)/n$ for this sample. 
iii) Repeat step ii) $m-1$ more times and obtain 
$\{\bar{y}^{B}_1,...,\bar{y}^{B}_m\mid \underline{y}\}\stackrel{iid}{\sim} H_{\underline{y}}.$
iv)
The $1-\alpha$ percentile bootstrap interval for $\mu$ based on $\bar{y}$ is
$C_{pN}=[\bar{y}_{(m_l)}^B, \bar{y}_{(m_u)}^B].$ 

We take a close look at $H_{\underline{y}}$. 
First, $H_{\underline{y}}$ depends on 
$\underline{y}$, not the complete and sufficient statistic $\bar{y}$ under the normal model. e.g., when $\underline{y}=\{0,1\}$ and $\underline{y}'=\{0.1,0.9\}$, $\bar{y}=\bar{y'}$  but   $H_{\underline{y}}\not = H_{\underline{y}'}$. 
Secondly, without loss of generality, assume, for two points
$\underline{y}^B=\{y_1^B,...,y_n^B\}$  and 
$\underline{y}^{'B}=\{y_1^{'B},...,y_n^{'B}\}$ in $S_{\underline{y}}$,  
\begin{equation}
\label{permutation}
\bar{y}^B 
=\bar{y'}^{B} \,\ \mbox{if and only if $\underline{y'}^{B}$ is a permutation of $\underline{y}^{B}$.}
\end{equation}
This indeed occurs with probability one. 
Let $S(n,t)$ be the number of terms in polynomial $(\sum_{i=1}^n x_i)^t$, which can be determined inductively following  $S(n,t)=\sum_{j=0}^t S(n-1, j)$,  $S(1,n)=1$ and $S(n,0)=1$.
Then, the estimator $\bar{y}^B$ assumes $S(n,n)$ distinct values over $S_{\underline{y}}$ due to (\ref{permutation}). For examples, $S(2,2)=3$ and $S(3,3)=10$. For a point $\underline{y}^B=\{y_1^B,...,y_n^B\}$ in $S_{\underline{y}}$, let $k_j$ be the number of elements in set $\{y_i^B: y_i^B=y_j\}$ for $j=1,...,n$. Then, $\sum_{j=1}^n k_j=n$ and the conditional PMF and CDF of the estimator  $\bar{y}^B\mid \underline{y}$ are 
\begin{equation}\label{probability mass function-cumulative distribution functions-mean}
p_{\underline{y}}(\bar{y}^B)={n^{-n} {n \choose k_1 ... k_n} } 
\,\ \mbox{and}\,\
H_{\underline{y}}(x)=n^{-n}\sum_{\bar{y}^B\leq x} { {n \choose k_1 ... k_n} }, 
\end{equation}
respectively.
The details are illustrated in Example A1 in the Appendix for $n=3$.

The coverage probability for $C_{pN}$ is equal to 
\begin{equation}\label{cp-c_pN}
\begin{split}
&Cover_{C_{pN}}(\mu)
=Cover_{C_{pN}}(0)\\
&= \int...\int [F_B(m_u-1,m,H_{\underline{y}}(0))  - F_B(m_l-1,m,H_{\underline{y}}(0))] \Pi_{i=1}^n \{ p_N(y_i,0,\sigma)dy_i\}.
\end{split}		
\end{equation}
It is easy to see that $H_{\underline{y}}(0)=1$ if $y_{(n)}<0$ and $H_{\underline{y}}(0)=0$ if $y_{(1)}>0$. So, for any $\alpha\in (0,1)$,
\begin{equation}\label{cp-upp}
Cover_{C_{pN}}(\mu)\leq 1-P(y_{(n)}<0\cup y_{(1)}>0)=1- 2^{-(n-1)}.
\end{equation}
Hence, the confidence coefficient cannot be $1-\alpha$ for a large $1-\alpha$. e.g., the  $C_{pN}$ with a $95\%$ confidence coefficient does not exist for $n=5$.
The expected length of $C_{pN}$ is equal to
\begin{equation}\label{el-c_pN}
\begin{split}
EL_{C_{pN}}(\mu) &=\int_{-\infty}^{+\infty}...\int_{-\infty}^{+\infty} \{\sum_{s=1}^{K-1}(a_{s+1}-a_s)[F_B(m_u-1,m,H_{\underline{y}}(a_s))\\ & \hspace{1.8in} -F_B(m_l-1,m,H_{\underline{y}}(a_s)) ]\}\{\Pi_{i=1}^n p_N(y_i,0,\sigma) dy_i\},
\end{split}
\end{equation}
where,  for a given $\underline{y}$,  $\{a_s\}_{s=1}^K$ are the possible values of $\bar{y}_j^B$ from the smallest to the largest with $K=S(n,n)$. The proof is in the Appendix.

For illustration, we derive $Cover_{C_{pN}}$ and $EL_{C_{pN}}$ for $n=2$ in Example A2 in the Appendix.
For a general $n$, it is difficult to derive them. Instead, we use the simulations to compute them.
Fortunately, since $Cover_{C_{pN}}(\mu)$ and $EL_{C_{pN}}(\mu)$ do not depend on $\mu$, simulating them at a single point of $\mu=0$ is enough. So, the simulation study here is trusted and complete because it covers all parameter configurations.  This, however, does not occur in most of the simulation studies. 

Table~\ref{tab-normal} contains the simulated coverage probability and expected length of $1-\alpha$ interval $C_{pN}$ using 100000 replications for each case. The three concerns on the interval $C_N$ remain for $C_{pN}$ but become worse. In addition, the nonparametric interval $C_{pN}$ is narrower than the parametric version $C_N$ and the $z$-interval under the normal model -- another illogical conclusion.


\subsection{The parametric bootstrap intervals of a normal mean $\mu$ based on the sample median}

A normal mean can also be estimated by $Med$,  the sample median.
For simplicity, assume $n$ is odd. Let $a=[n/2]+1$. Then $Med=y_{(a)}$.
The parametric bootstrap interval for $\mu$ based on $Med$ is obtained by applying Theorem~\ref{thm-1} to $Q=Med$ and is equal to 
$C_{NM}=[Med_{(m_l)}^B, Med_{(m_u)}^B].$
The coverage probability for $C_{NM}$ is the constant in (\ref{cp-normal-q}). Let
$F_{Be}(x,a,b)$ be the CDF of Beta distribution with parameters $a$ and $b$. The CDF of 
$Med$ for $\mu=0$ and $\sigma=1$ is
$F(x,0,1)=F_{Be}(F_N(x,0,1),a,n-a+1).$
Thus,
the expected length of $C_{NM}$, following (\ref{el-normal-q}), is given by
\begin{eqnarray*}
	EL_{C_{NM}}(\mu)
	=  m \sigma \int_0^1 F_N^{-1}(F_{Be}^{-1}(z,a,n-a+1),0,1) [
	p_B(m_u-1,m-1,z) 
	-p_B(m_l-1,m-1,z)]dz.
\end{eqnarray*}

Table 2 contains the coverage probability and expected length for 
$C_{NM}$. The concerns a) and b) for $C_N$ are also true for $C_{NM}$. i.e., its coverage probability is less than the nominal level and is not increasing in $m$. As expected, $C_{NM}$ is much wider than $C_N$ since the latter is based on the UMVUE $\bar{X}$ for $\mu$.  
In the next section, we derive the percentile bootstrap interval for $\mu$ based on $Med$ in a general setting.

\section{The percentile bootstrap interval for the median of a symmetric and continuous population based on $MED$}\,\

A random observation $Y$ has a symmetric PDF $f(y)$ and a CDF $F(y)$, where $f$ belongs to a parametric or nonparametric distribution family ${\cal F}$. Suppose the median $\mu(f)$ is the unique 50-th percentile for $Y$. So,  
$\mu(f)$ is the unique solution for $\int _{-\infty}^{\mu} f(y)dy=0.5$. 

Derive the percentile bootstrap interval for $\mu(f)$ as follows.
i) Generate an i.i.d. sample of size $n$ (odd): $\underline{y}=\{y_1,...,y_n\}\sim f(y)$. 
Without loss of generality, we assume all $y_i$'s are distinct because it occurs with probability one.
ii) For this $\underline{y}$, consider $S_{\underline{y}}$ in (\ref{boot-s}). Define the random variable $Med^B=y_{(a)}^B$ on $S_{\underline{y}}$ with a CDF
$H_M$. Generate an i.i.d. bootstrap sample $\{y_{11}^B,...,y_{n1}^B\}$ from $\underline{y}$ 
and compute $Med_1^B$. 
iii) Repeat ii) for $m-1$ more times and obtain 
$	\{Med^{B}_1,...,Med^{B}_m\mid \underline{y}\}\stackrel{iid}{\sim} H_M.$
iv)
The $1-\alpha$ percentile bootstrap interval for $\mu(f)$ based on $Med$  is
$C_{pM}=[Med_{(m_l)}^B, Med_{(m_u)}^B].$

We show in the Appendix that the range of $Med^B$ is $\{y_{(i)}\}_{i=1}^n$ and 
\begin{equation}
\label{cumulative distribution functions-med}
H_M(y_{(i)})=\sum_{w=0}^{a-1} {n\choose w} ({n-i\over n})^w ({i\over n})^{n-w}.
\end{equation}

\begin{thm}\label{thm-median} Under the conditions in the first paragraph of this section, the coverage probability and expected length for $C_{pM}$, as functions of $f\in {\cal F}$, 
	are given below
	\begin{equation}\label{cp-c_pM}
	\begin{split}
	Cover_{C_{pM}}(f) 
	&= \sum_{i=0}^n  [F_B(m_u-1,m,F_B(a-1,n,{n-i\over n})) \\& \hspace{0.5in} -F_B(m_l-1,m,F_B(a-1,n,{n-i\over n}))] p_B(i,n,0.5);
	\end{split}		
	\end{equation}
	\begin{equation}\label{el-c_pM}
	\begin{split}
	&			EL_{C_{pM}}(f) 
	=n \sum_{i=1}^{n-1}[ F_B(m_u-1,m,F_B(a-1,n,{n-i\over n}))   -F_B(m_l-1,m,\\ & \hspace{0.9in} F_B(a-1,n,{n-i\over n}))]
	\int_0^1 F^{-1}(z)[p_B(i,n-1,z)-p_B(i-1,n-1,z)]dz.
	\end{split}
	\end{equation}
\end{thm}

The coverage probability of $C_{pM}$ is constant in $f$. Also, 
the same conclusion as (\ref{cp-upp}) holds for interval $C_{pM}$.  
When ${\cal F}$ is the normal family $N(\mu,\sigma^2)$ with a known (or unknown) $\sigma^2$,  
$C_{pM}$ can be used to estimate $\mu$. The coverage probability is given in (\ref{cp-c_pM}) with $f$ replaced by $\mu$ and the expected length changes to
\begin{equation}\label{el-c_pNM}
\begin{split}
&			EL_{C_{pM}}(\mu) 
=n\sigma \sum_{i=1}^{n-1}[F_B(m_u-1,m,F_B(a-1,n,{n-i\over n}))   -F_B(m_l-1,m,\\ & \hspace{0.8in} F_B(a-1,n,{n-i\over n}))]
\int_0^1 F^{-1}_N(z,0,1)[p_B(i,n-1,z)-p_B(i-1,n-1,z)]dz.
\end{split}
\end{equation}

 Table~\ref{tab-normal} contains the coverage probability and expected length for $C_{pM}$. The comparison between $C_{pM}$ and its parametric version $C_{NM}$ is not as consistent as the one between $C_{pN}$ and $C_N$. 
 In general, $C_{pM}$ has a smaller expected length than $C_{NM}$ when $n$ is not small.
 This again questions the usage of a parametric model for statistical inferences.

\section{The parametric bootstrap interval of a normal mean based on the sample mean when $\sigma^2$ is unknown}

An i.i.d. sample $\underline{y}=\{y_1,...,y_n\}$ is observed from $N(\mu,\sigma^2)$.
The parametric bootstrap interval for $\mu$ based on $\bar{y}$ is derived as follows: 
i) Generate an i.i.d. sample: $\{y_1,...,y_n\} \sim N(\mu,\sigma^2)$. 
Compute $\bar{y}= \sum_{i=1}^{n}x_i/n$ and $s^2=\sum_{i=1}^n(y_i-\bar{y})^2/n$, the MLEs of $\mu$ and $\sigma^2$. 
ii) For this pair $(\bar{y},s^2)$, generate an i.i.d. sample
$		\{y_{11}^B,...,y_{n1}^B\mid \bar{y},s^2\}\sim N(\bar{y},s^2).$
Compute $\bar{y}^{B}_1=\sum_{i=1}^{n}y_{i1}^B/n$. Then $\bar{y}^{B}_1\mid \bar{y},s^2\sim N(\bar{y},{s^2/ n}).$
iii) Repeat ii) for $m-1$ more times and obtain an i.i.d. sample
$	\{\bar{y}^{B}_1,...,\bar{y}^{B}_m\mid \bar{y},s^2\} \sim N(\bar{y},{s^2/ n}).$
iv) The $1-\alpha$ parametric bootstrap interval for $\mu$ is
$C_{Nu}=[\bar{y}_{(m_l)}^B, \bar{y}_{(m_u)}^B].$

As shown in the Appendix, the coverage probability for $C_{Nu}$ is given by
\begin{equation}\label{cp-C_Nu}
\begin{split}
Cover_{C_{Nu}}(\mu,\sigma)
&= E[F_B(m_u-1,m,F_N(- T\sqrt{{n\over n-1}} ,0,1))\\ &\hspace{0.5in}
-F_B(m_l-1,m,F_N(-T\sqrt{{n\over n-1}},0,1))],
\end{split}
\end{equation}
where $T$ follows the t-distribution with $n-1$ degrees of freedom. So, this function is independent of $(\mu,\sigma)$ but depends on $(1-\alpha, n, m)$.
The expected length of $C_{Nu}$ is equal to
\begin{equation}\label{el-C_Nu}
EL_{C_{Nu}}(\mu,\sigma)
= {\sqrt{2}\over\sqrt{n}}{ \Gamma({n\over 2})\over \Gamma({n-1\over 2}) }A(n,m,\alpha)\sigma\stackrel{def}{=}B(n)A(n,m,\alpha)\sigma, 
\end{equation}			
where $\Gamma(\cdot)$ is the gamma function and $A(n,m,\alpha)$ is given in (\ref{elnormal}).

Similar to the domination of $C_N$ over the $z$-interval in length, Table~\ref{tab-normal} shows that $C_{Nu}$ also dominates the $t$-interval, an illogical conclusion. Furthermore, 
$B(n)$, the ratio of $EL_{C_{Nu}}(\mu,\sigma)/EL_{C_N}(\mu)$, is always less than 1, indicating that $C_{Nu}$ is narrower than $C_N$ for the same nominal level $1-\alpha$. 
This is counter-intuitive: $C_{Nu}$,  involving two unknown parameters $\mu$ and $\sigma^2$, should be wider than $C_N$, involving only one unknown parameter $\mu$, but we have the opposite. e.g., when $(1-\alpha, n,m)=(0.9,31,100)$, $C_N$ has 
an expected length 0.5716, while  $C_{Nu}$ has 
an expected length 0.5576. On the contrary, the $t$-interval is always wider than the $z$-interval. In Table~\ref{tab-normal} we find, in the order of the expected length,  
\begin{equation}
\label{three}
C_{pN}\leq C_{N}\leq \mbox{the $z$-interval}\,\ \mbox{and}\,\ 
C_{Nu}\leq C_{N}\leq \mbox{the $z$-interval}
\end{equation}
even when the underlying distribution is normal. The above relationships are due to the same nominal level for three intervals. 


\section{Discussion}						

It is a common sense that a meaningful comparison must have the same objective baseline. In practice, the nominal level $1-\alpha$ is generally used to construct confidence intervals, then a comparison in expected length is conducted among intervals of the same nominal level. Whether the nominal level is a good baseline becomes an important issue.  

The usage of nominal level is simple but, as shown in our paper, causes 
many problems for bootstrap intervals summarized below: i) the coverage probability may not increase in the nominal level (e.g., $C_{wi}$); ii) a larger 
$m$,  meaning more efforts and more information, may not yield a shorter interval (e.g., $C_{Nu}$ and $C_{NM}$); iii) the choice of $m=50$ for ``fairly good standard error estimates'' is not acceptable; iv) the expected length does not increase in the number of unknown parameters (e.g., $C_{N}$ vs $C_{Nu}$); v) the inference based on a parametric model is worse than that based on a nonparametric model (e.g., $C_N$ vs $C_{pN}$,  $C_{NM}$ vs $C_{pM}$); vi) the unacceptable relationship (\ref{three}): the parametric bootstrap intervals $C_N$ and $C_{Nu}$ are surprisingly narrower than the optimal $z$-interval; vii)
the confidence coefficient may be zero for any nominal level (e.g., $C_{wa}$, $C_{wi}$, $C_{d}$ and $C_{\theta}$). 

In practice, the nominal level and the confidence coefficient 
 are often treated the same but are truly different as shown in Tables \ref{tab-wilson} and \ref{tab-normal}. The former is subjective because it is predetermined. A better choice 
in interval comparison is to compare the length among those intervals with the same confidence coefficient or the same area under coverage probability curve (not the same nominal level).  This will solve the problems listed above at least in a certain degree. Figure~\ref{fig4} gives an example of using the area under curve for comparison, another example, Example A3, of using the confidence coefficient is given in the Appendix.

With the proposed general method in Section 2, we are able to calculate coverage probability and expected length analytically, then study the finite-sample properties of bootstrap intervals without the vague assumption of ``large samples''. The bootstrap interval is now being examined under the newly invented ``microscope''. In particular, it is possible to conduct
a fair comparison of bootstrap intervals using an objective baseline. 
Due to the poor performance of all simple bootstrap intervals discussed in the paper, it is hard to expect any good performance of the bootstrap interval in a complicated case.
The conclusion drawn in any of its application must be carefully evaluated and interpreted.

\section{APPENDIX}
{\it Proof of Equation (\ref{el-wild}).}
	Following (\ref{distphat}) and (\ref{order-CDF-p-g}), the conditional CDF of the $j$-th order statistic $n\hat{p}^B_{(j)}\mid y$ is  
	$F_j(x)=1-F_B(j-1,m,F_B(x,n,{y/n}))\,\ \mbox{for any} \,\ x\in R^1.$
	Also, $n\hat{p}^B_j$ assumes values, $0,...,n.$ So does $n\hat{p}^B_{(j)}$. Then,
	\begin{eqnarray*}
		E(\hat{p}^B_{(j)}\mid y)&=&
		{1\over n}\sum_{x=0}^{n} x(F_j(x)-F_j(x^-))
		={1\over n}\sum_{x=0}^{n} x(F_j(x)-F_j(x-1))\\
		&=&1-{1\over n}\sum_{x=0}^{n-1} F_j(x)
		= {1\over n}\sum_{x=0}^{n-1}F_B(j-1,m,F_B(x,n,{y\over n})).
	\end{eqnarray*}
	Therefore, (\ref{el-wild}) follows
	$ EL_{C_{wa}}(p)=\sum_{y=0}^n \{E(\hat{p}^B_{(m_u)}\mid y)-E(\hat{p}^B_{(m_l)}\mid y)\}p_B(y,n,p). $~\raisebox{.5ex}{\fbox{}}

{\it Proof of Equation (\ref{ecd}).}
	\begin{eqnarray*}
		E(\hat{d}^B_{(j)}\mid x,y)&=& \sum_{s=1}^K d_s(F_j(d_s)-F_j(d_s^-))
		=\sum_{s=2}^K d_s(F_j(d_s)-F_j(d_{s-1}))+d_1F_j(d_1)\\
		&=& \sum_{s=1}^K d_s F_j(d_s)- \sum_{s=2}^K d_sF_j(d_{s-1})	= \sum_{s=1}^K d_s F_j(d_s)- \sum_{s=1}^{K-1} d_{s+1}F_j(d_s)\\
		&=&d_K-\sum_{s=1}^{K-1}(d_{s+1}-d_s)F_j(d_s)
		=d_K-\sum_{s=1}^{K-1}(d_{s+1}-d_s)(1-F_B(j-1,m,H_{x,y}^D(d_s)))\\
		&=& d_1+ \sum_{s=1}^{K-1}(d_{s+1}-d_s)F_B(j-1,m,H_{x,y}^D(d_s)). ~\raisebox{.5ex}{\fbox{}}
	\end{eqnarray*}

{\it Proof of Theorem~\ref{thm-1}.}
	\begin{eqnarray*}
		Cover_{C_Q}(\mu)
		&=& E_{Q \sim  p(q,\mu,\sigma)\}}[{\rm pr}(Q_{(m_l)}^B-\mu\leq 0 \leq Q_{(m_u)}^B-\mu\mid Q )]
		\\&=& E_{\{Q \sim p(q,0,\sigma)\}}[{\rm pr}(Q_{(m_l)}^B \leq 0 \leq Q_{(m_u)}^B\mid Q)]=	Cover_{C_Q}(0).
	\end{eqnarray*}
	Let $F_j$ be the conditional CDF of order statistic $Q^B_{(j)}\mid Q$. 
	Note
	$$F(0,Q,\sigma)
	=1-F(Q,0,\sigma) \,\ \mbox{and}\,\
	F_j(0)
	=1-F_B(j-1,m,1-F(Q,0,\sigma)).$$
	Therefore,
	\begin{eqnarray*}
		\begin{split}
			&Cover_{C_Q}(0) 
			=\int_{-\infty}^{+\infty}[F_B(m_u-1,m,1-F(q,0,\sigma))
			-F_B(m_l-1,m,1-F(q,0,\sigma))] p(q,0,\sigma) dq
			\\
			&\hspace{0.75in}=\int_{0}^{1}[F_B(m_u-1,m,1-z)-F_B(m_l-1,m,1-z)]  dz\\ 
			&\hspace{0.75in}= \sum_{k=m_l}^{m_u-1}\int_0^1 p_B(k,m,1-z)dz=\sum_{k=m_l}^{m_u-1}{1\over m+1}={m_u-m_l\over m+1}.
		\end{split}
	\end{eqnarray*}
	
	The expected length of $C_Q$ is equal to
	\begin{eqnarray*}
		EL_{C_Q}(\mu)&=&
		E_{\{Q \sim  p(q,\mu,\sigma)\}}[E(Q^B_{(m_u)}-Q^B_{(m_l)}\mid Q)]=E_{\{Q \sim p(q,0,
			\sigma) \}}[E\{Q^B_{(m_u)}-Q^B_{(m_l)}\mid Q\}]\\ &=&\int_{-\infty}^{+\infty} (E(Q^B_{(m_u)}\mid q)-E(Q^B_{(m_l)}\mid q))p(q,0,\sigma) dq,
	\end{eqnarray*}
	which is independent of $\mu$.
	Note that the conditional CDF of $Q^B_j\mid Q$ is $F(x,Q,\sigma)$. Following (\ref{order-CDF-p-g}), the conditional PDF of order statistic $Q^B_{(j)}\mid Q$ is 
	$$f_j(x)=m \,\ p(x,Q,\sigma) p_B(j-1,m-1,F(x,Q,\sigma)) \,\ \mbox{for any} \,\ x\in R^1.$$
	Let $z=F(x,Q,\sigma)$. Then, $x=F^{-1}(z,0,1)\sigma+Q$ due to (\ref{f-form}). Therefore, (\ref{el-normal-q}) follows
	\begin{eqnarray*}
		EL_{C_Q}(\mu)
		&=& \int_{-\infty}^{+\infty} [\int_{0}^{1} m(F^{-1}(z,0,1)\sigma+q) (
		p_B(m_u-1,m-1,z) \\ && \hspace{2.1in}
		-p_B(m_l-1,m-1,z))dz]p(q,0,\sigma) dq\\
		&=& m\sigma \int_{0}^{1} F^{-1}(z,0,1) (
		p_B(m_u-1,m-1,z) -p_B(m_l-1,m-1,z))dz. ~\raisebox{.5ex}{\fbox{}}
	\end{eqnarray*}

{\it Example A1.}\label{ex-1}
	When $n=3$ and $\underline{y}=(0,1,5)$, $S_{\underline{y}}$ and $\bar{y}^B$ are given in Table~\ref{tab1}. 	~\raisebox{.5ex}{\fbox{}}
\begin{table}[h]
	\begin{center}
	\caption{$S_{\underline{y}}$ (with 27 points), $\bar{y}^B$ (with 10 values), $(k_1,k_2,k_3)$, the PMF and CDF of $\bar{y}^B$. \label{tab1}}{
	\begin{tabular}{ccccc}\hline
		$\underline{y}^B$ & $\bar{y}^B$ & $(k_1,k_2,k_3)$ &  $p_{\underline{y}}(\bar{y}^B)$ in (\ref{probability mass function-cumulative distribution functions-mean})&  $H_{\underline{y}}(\bar{y}^B)$ in (\ref{probability mass function-cumulative distribution functions-mean})\\ 
		\hline
		$(0,0,0)$                    & 0            & (3,0,0)& $1/27$ & $1/27$\\
		$ (0,0,1), (0,1,0), (1,0,0)$ & $1/3$ & (2,1,0)& $3/27$ & $4/27$\\
		$ (0,1,1), (1,0,1), (1,1,0)$ & $2/3$ & (1,2,0)& $3/27$ & $7/27$\\
		$ (1,1,1)                  $ & $1$          & (0,3,0)& $1/27$ & $8/27$\\
		$ (0,0,5), (0,5,0), (5,0,0)$ & $5/3$ & (1,0,2)& $3/27$ & $11/27$\\
		$ (0,1,5), (0,5,1), (1,0,5)$ & $2$          & (1,1,1)& $6/27$ & $17/27$\\
		$ (1,5,0), (5,0,1), (5,1,0)$ &              &        &               & \\
		$ (1,1,5), (1,5,1), (5,1,1)$ & $7/3$ & (0,2,1)& $3/27$ & $20/27$\\
		$ (0,5,5), (5,0,5), (5,5,0)$ & $10/3$ & (1,0,2)& $3/27$ & $23/27$\\
		$ (1,5,5), (5,1,5), (5,5,1)$ & $11/3$ & (0,1,2)& $3/27$ & $26/27$\\
		$ (5,5,5)                 $ & $5$           & (0,0,3)& $1/27$ & $1$\\ \hline
\end{tabular}}
	\end{center}
\end{table}


{\it Example A2.}\label{ex-2}
	We compute $Cover_{C_{pN}}$ and $EL_{C_{pN}}$ when $n=2$ and $\underline{y}=\{y_1,y_2\}$. Then, 
	$$H_{\underline{y}}(0)
	= \sum_{ \bar{y}^B\leq 0} {
		{2\choose k_1}}/4=\left\{ \begin{tabular}{ll}
	$0$& if $(y_1,y_2)\in A_1=\{0<y_{(1)}\}$;\\
	$0.25$ & if $(y_1,y_2)\in A_2=\{y_{(1)}\leq 0<{y_{(1)}+y_{(2)}\over 2}\}$;\\
	${0.75}$ & if $(y_1,y_2)\in A_3=\{ {y_{(1)}+y_{(2)}\over 2} \leq 0< y_{(2)}\}$;\\
	$1$          & if $(y_1,y_2)\in A_4=\{0\geq y_{(2)}\}$.
	\end{tabular}
	\right. 
	$$
	Following (\ref{cp-c_pN}), since the integrand below is zero on $A_1$ and $A_4$, 
	\begin{eqnarray*}
		Cover_{C_{pN}}(\mu) &=&  \int\int_{A_1\cup A_2\cup A_3\cup A_4} [F_B(m_u-1,m,H_{\underline{y}}(0))- F_B(m_l-1,m,H_{\underline{y}}(0)) ] \Pi_{i=1}^2p_N(y_i,0,\sigma)dy_i	\\
		&=& 2^{-1}[F_B(m_u-1,m,0.25)-F_B(m_l-1,m,0.25)].
	\end{eqnarray*}
	\\
	The true confidence coefficient of any $1-\alpha$ interval $C_{pN}$ is no larger than 0.5 due to (\ref{cp-upp}) even if the nominal level $1-\alpha$ is set to be $0.9, 0.95$ or 0.99. 		
	Each $\bar{y}^B_j$ assumes $3(=K)$ values: $a_1=y_{(1)}$, $a_2=((y_{(1)}+y_{(2)})/2$ and $a_3=y_{(2)}$, and $H_{\underline{y}}(a_1)=0.25$, $H_{\underline{y}}(a_2)=0.75$ and $H_{\underline{y}}(a_3)=1$. Following (\ref{el-c_pN}),
	$$
	EL_{C_{pN}}(\mu) =2\sigma\{F_B(m_u-1,m,0.25) -F_B(m_l-1,m,0.25) \} \int_0^1 F_N^{-1}(z,0,1)(2z-1)dz. 	~\raisebox{.5ex}{\fbox{}}$$

{\it Proof of Equation (\ref{el-c_pN}).}
	\begin{eqnarray*}
		EL_{C_{pN}}(\mu)
		&=&\int_{-\infty}^{+\infty}...\int_{-\infty}^{+\infty} [E(\bar{y}^B_{(m_u)}\mid \underline{y})-E(\bar{y}^B_{(m_l)}\mid \underline{y})]\Pi_{i=1}^n p_N(y_i,0,\sigma) dy_i.
	\end{eqnarray*}
	For the observed $\underline{y}$, the conditional CDF of $\bar{y}^B_{(j)}\mid \underline{y}$  for a given $j\in [1,m]$ is
	$$F_j(x)=1-F_B(j-1,m,H_{\underline{y}}(x))\,\ \mbox{for any} \,\ x\in R^1.$$ Then, (\ref{el-c_pN}) follows
	\begin{eqnarray*}
		E(\bar{y}^B_{(j)}\mid \underline{y})
		&=& a_1+ \sum_{s=1}^{K-1}(a_{s+1}-a_s)F_B(j-1,m,H_{\underline{y}}(a_s)). ~\raisebox{.5ex}{\fbox{}}
	\end{eqnarray*}

{\it Proof of Equation (\ref{cumulative distribution functions-med}).}
	Since the sample size $n$ is odd, the range of $Med^B$ is $\{y_{(i)}\}_{i=1}^n$. Note
	$ H_M(y_{(i)})= {\rm pr}(Med^B\leq y_{(i)})={\rm pr}(w_i<a),$
	where $w_i$ is the number of $y_j^B$'s in a sequence of $\{y_1^B,...,y_n^B\}$ larger than $y_{(i)}$. For each $y_j^B$ in the sequence, $y_j^B$ is either larger than $y_{(i)}$ (success), or no larger than $y_{(i)}$ (failure), ${\rm pr}(y_j^B>y_{(i)})={(n-i) / n}$ due to the distinct $y_i$'s; this probability remains unchanged in $j$ and all $y_j^B$'a are independent. Thus,
	$w_i$ follows $Bino(n,{(n-i)/ n})$ and (\ref{cumulative distribution functions-med}) is true.	
~\raisebox{.5ex}{\fbox{}}

{\it Proof of Theorem~\ref{thm-median}.}
	Note
	$ H_M(\mu(f))
	=F_B(a-1,n, {(n-i)/n}),$
	and, for a given $\underline{y}=\{y_1,...,y_n\}$ with distinct $y_i$'s, the integer $i$ is equal to the number of $y_j$'s less than or equal to $\mu(f)$. 
	Let $R_n=\{(y_1,...,y_n)\in R^n: \mbox{$y_j$'s are distinct}\}$. Then, $R_n$ has a partition below:
	$$R_n=\{\mu(f)<y_{(1)}\} \cup \{y_{(1)}\leq \mu(f)<y_{(2)}\}\cup ...\cup \{y_{(n)\leq \mu(f)}\}\stackrel{def}{=} \cup_{i=0}^n R_{n,i}.$$ Also, $i$ follows $Bino(n,0.5)$ because $y_j$ has a symmetric $f(y)$; on each $R_{n,i}$, $$H_M(\mu(f))=F_B(a-1,n,{n-i\over n}), \int...\int_{R_{n,i}} \Pi_{u=1}^n f(y_u)dy_u=p_B(i,n,0.5).$$
	Therefore,	\begin{eqnarray*}
		&&Cover_{C_{pM}}(f)
		=E_{\{y_i's\stackrel{iid}{\sim} f(y)\}}[P_{\underline{y}}(Med_{(m_l)}^B\leq \mu(f) \leq Med_{(m_u)}^B\mid \underline{y} )]\\ 
		&=& \int...\int [F_B(m_u-1,m,H_M(\mu(f)^-))-	F_B(m_l-1,m,H_M(\mu(f))) ] \Pi_{u=1}^n f(y_u) dy_u\\ 
		&=& \sum_{i=0}^n \int...\int_{R_{n,i}} [F_B(m_u-1,m,H_M(\mu(f)))- F_B(m_l-1,m,H_M(\mu(f))) ] \Pi_{u=1}^n f(y_u)dy_u\\
		&=& \sum_{i=0}^n   [F_B (m_u-1,m,F_B (a-1,n,{n-i\over n} ) )- F_B (m_l-1,m,F_B (a-1,n,{n-i\over n} ) )  ] p_B(i,n,0.5). 
	\end{eqnarray*}
	
	For an observed $\underline{y}$, the CDF of $Med^B_{(j)}\mid \underline{y}$ is
	$$F_j(x)=1-F_B(j-1,m,H_M(x))\,\ \mbox{for any} \,\ x\in R^1.$$
	Note that $Med^B$ assumes values $\{y_{(i)}\}_{i=1}^n$.  Then,
	\begin{eqnarray*}
		E(Med^B_{(j)}\mid \underline{y})
		&=& y_{(1)}+ \sum_{i=1}^{n-1}(y_{(i+1)}-y_{(i)})F_B(j-1,m,H_M(y_{(i)})).
	\end{eqnarray*}
	For any fixed $i$, $H_M(Y_{(i)})$ is a constant given in (\ref{cumulative distribution functions-med}) and
	\begin{eqnarray*}
		E(Y_{(i)})=n \int_{R^1} x f(x)p_B(i-1,n-1,F(x))dx
		= n \int_0^1 F^{-1}(z)p_B(i-1,n-1,z)dz.
	\end{eqnarray*}
	Therefore,
	\begin{eqnarray*}
		&&EL_{C_{pM}}(f)=
		E_{\underline{y}
			\stackrel{iid}{\sim} f(y)\}}E\{(Med^B_{(m_u)}-Med^B_{(m_l)})\mid \underline{y}\}\\ 
		&=&		\sum_{i=1}^{n-1}E [(Y_{(i+1)}-Y_{(i)})\{F_B(m_u-1,m,H_M(Y_{(i)})) -F_B(m_l-1,m,H_M(Y_{(i)}))\}]\\
		&=&\sum_{i=1}^{n-1}[F_B(m_u-1,m,H_M(Y_{(i)})) -F_B(m_l-1,m,H_M(Y_{(i)}))] (EY_{(i+1)}-EY_{(i)}), 
	\end{eqnarray*}	
	and	(\ref{el-c_pM}) is	established .   
~\raisebox{.5ex}{\fbox{}}

{\it Proof of Equations~(\ref{cp-C_Nu}) and (\ref{el-C_Nu}).}
	Let $F^*$ and $F_j$ be the conditional CDFs of $\bar{y}^B_j$ and the order statistic $\bar{y}^B_{(j)}$ for given $(\bar{y},s^2)$, respectively. So,
	$$F^*(0)
	=F_N(0,\bar{y},\sqrt{{s^2\over n}})
	=F_N(-{\bar{y}\over \sqrt{{s^2\over n}} },0,1)
	\,\ \mbox{and}\,\
	F_j(0)
	=1-F_B(j-1,m,F^*(0)).
	$$
	The coverage probability for $C_{Nu}$ satisfies
	\begin{eqnarray*}
		&&Cover_{C_{Nu}}(\mu,\sigma)=E_{\{\bar{y}\sim N(\mu,{\sigma^2\over n}),s^2\sim {1\over n}\sigma^2 \chi^2_{n-1} \}}[P(\bar{y}_{(m_l)}^B\leq \mu \leq \bar{y}_{(m_u)}^B\mid \bar{y},s^2 )]\\
		&=& E_{\{\bar{y}\sim N(0,{\sigma^2\over n}),s^2\sim {1\over n}\sigma^2 \chi^2_{n-1}\}}[P(\bar{y}_{(m_l)}^B \leq 0 \leq \bar{y}_{(m_u)}^B\mid \bar{y},s^2)]\\
		&=&  E_{\bar{y}\sim N(0,{\sigma^2\over n}),s^2\sim {1\over n}\sigma^2 \chi^2_{n-1}\}}[P(\bar{y}_{(m_l)}^B \leq 0\mid \bar{y},s^2)- P( \bar{y}_{(m_u)}^B\leq 0\mid \bar{y},s^2)]\\
		&=& E[F_B(m_u-1,m,F_N(-{\bar{y}\over \sqrt{{s^2\over n}} },0,1))-F_B(m_l-1,m,F_N(-{\bar{y}\over \sqrt{{s^2\over n}} },0,1))].
	\end{eqnarray*}
	Then, (\ref{cp-C_Nu}) follows ${\bar{y}/\sqrt{{(s^2/ n)}}}\sim T\sqrt{{(n/(n-1))}}$.
	
	Note that
	$$	EL_{C_{Nu}}(\mu,\sigma)=
	E_{\{\bar{y}\sim N(0,{\sigma^2\over n}),s^2\sim {1\over n}\sigma^2 \chi^2_{n-1}\}}E[\bar{y}^B_{(m_u)}-\bar{y}^B_{(m_l)}\mid \bar{y},s^2],$$
	independent of $\mu$,
	and
	the conditional PDF of $\bar{y}^B_{(j)}\mid \bar{y}, s^2$ is
	$$f_j(x)=m \,\  p_N(x,\bar{y},\sqrt{{s^2\over n}}) p_B(j-1,m-1,F_N(x,\bar{y},\sqrt{{s^2\over n}}))\,\ \mbox{for any} \,\ x\in R^1.$$
	Let $z=F_N(x,\bar{y},\sqrt{{(s^2/ n)}})$. 
	Then, $x=F_N^{-1}(z,0,1){s/ \sqrt{n}}+\bar{y}$. Therefore, 
	\begin{eqnarray*}
		&&E[\bar{y}^B_{(m_u)}-\bar{y}^B_{(m_l)}\mid \bar{y},s^2] 
		=  \int_{-\infty}^{+\infty} x(f_{m_u}(x)-f_{m_l}(x))dx \\
		&=& \int_{-\infty}^{+\infty} m\,\ x \,\ p_N(x,\bar{y},\sqrt{{s^2\over n}})[ p_B(m_u-1,m-1,F_N(x,\bar{y},\sqrt{{s^2\over n}}))\\
		&&\ \ \ \ \ \ \ \ \ \ \ \ \ \ \ \ \ \ \ \ \ \ \ \ \ \ \ \ \ \ \ \ \ \ \ \ - p_B(m_l-1,m-1,F_N(x,\bar{y},\sqrt{{s^2\over n}}))]dx\\
		&=& \int_{0}^{1} m(F_N^{-1}(z,0,1){s\over \sqrt{n}}+\bar{y}) [
		p_B(m_u-1,m-1,z)-p_B(m_l-1,m-1,z)]dz\\
		&=& A(n,m,\alpha)s=A(n,m,\alpha){\chi_{n-1} \over \sqrt{n}} \sigma,
	\end{eqnarray*}
	where $A(n,m,\alpha)$ is given in (\ref{elnormal}), $\chi^2_k$ follows $\chi^2$-distribution with $k$ degrees of freedom. Then, (\ref{el-C_Nu}) is established by taking the expectation.
~\raisebox{.5ex}{\fbox{}}

{\it Example A3.} Table~\ref{tabCp_N} contains the comparison between the bootstrap interval $C_{pN}$ and the $z$-interval of level $1-\alpha$. The result for $C_{pN}$ is simulated based on 100000 replications for each case. At the first glance,  $C_{pN}$ seems better than the $z$-interval since the former is always shorter than the latter. However, the confidence coefficient of $C_{pN}$ is less than that of the $z$-interval. For a fair comparison, we further report the expected length of another $z$-interval, denoted by the z$^*$-interval, whose confidence coefficient is equal to that of the $1-\alpha$ interval $C_{pN}$. Then the conclusion flips over completely: the $z^*$-interval is narrower! Unfortunately, this misleading comparison between the
$1-\alpha$ interval $C_{pN}$ and the $1-\alpha$ z-interval is typically conducted.
The same occurs in other comparisons including $C_N$ vs the  $z$-interval and $C_{Nu}$ vs the $z$-interval. ~\raisebox{.5ex}{\fbox{}}

\begin{table}[H]
	\caption{The coverage probabilities $(CP)$ and expected lengths $(EL)$ for the $1-\alpha$ interval $C_{pN}$, the $1-\alpha$ $z$-interval and the $z^{\star}$-interval for different $(m, \alpha)$.}
	\begin{center}\vspace{-0.2in}
	\resizebox{0.8\textwidth}{0.2\textheight}
	{
		\begin{tabular}{l| cccl}
			&\multicolumn{3}{c}{($CP$ for $C_{pN}$, $EL$ for $C_{pN}$, $EL$ for $z^{\star}$-interval)} &$EL$ for $z$-interval \\
			$1-\alpha$ &$m=50$ & $m=100$ & $m=5000$ &    \\
			\hline
			 &\multicolumn{4}{c}{$n=5$}\\
			$0.99$   & (0.8598, 1.6114, 1.3193) & (0.8858, 1.7654, 1.4128) & (0.8984, 1.8316, 1.4643) & 2.3039 \\
			$0.95$   & (0.8083, 1.3711, 1.1677) & (0.8284, 1.4327, 1.2228) & (0.8372, 1.4512, 1.2484) & 1.7530 \\
			$0.9$    & (0.7730, 1.2239, 1.0806) & (0.7677, 1.2033, 1.0684) & (0.7874, 1.2489, 1.1149) & 1.4712 \\ 
			&\multicolumn{4}{c}{$n=30$}\\
			$0.99$   & (0.9508, 0.7991, 0.7182) & (0.9707, 0.5802, 0.7958) & (0.9826, 0.9121, 0.8684) & 0.9406 \\
			$0.95$   & (0.9073, 0.6601, 0.6139) & (0.9264, 0.6924, 0.6533) & (0.9354, 0.6966, 0.6748) & 0.7157 \\
			$0.9$    & (0.8674, 0.5802, 0.5492) & (0.8648, 0.5663, 0.5455) & (0.8830, 0.5854, 0.5724) & 0.6006 \\ 
			&\multicolumn{4}{c}{$n=50$}\\
			$0.99$   & (0.9555, 0.6258, 0.5683) & (0.9759, 0.6980, 0.6380) & (0.9865, 0.7146, 0.6987) & 0.7286 \\
			$0.95$   & (0.9127, 0.5163, 0.4836) & (0.9333, 0.5421, 0.5186) & (0.9414, 0.5453, 0.5349) & 0.5544 \\
			$0.9$    & (0.8746, 0.4536, 0.4335) & (0.8707, 0.4434, 0.4290) & (0.8898, 0.4579, 0.4518) & 0.4652  ~\raisebox{.5ex}{\fbox{}} \\
			
		\end{tabular}
	}\end{center}
\label{tabCp_N}
\end{table}

\vspace{-0.2in}
\setlength{\parindent}{-1.5em}

\begin{Large}
	{\bf References}
\end{Large}
\bigskip

Blyth, C. R. and Still, H. A. (1983), `` Binomial Confidence Intervals,'' \textit{Journal of the American Statistical Association}, 78, 108-116.
		
Casella, G. (1986), ``Refining Binomial Confidence Intervals,'' 
	\textit{Canadian  Journal of Statistics}, 14, 113-129.
		
Casella, G., and Berger, R. L. (2002), \textit{``
		Statistical Inference,}'' 2nd ed. Duxbury Press, Pacific Grove, CA.
	
Clopper, C. J., and Pearson, E. S. (1934). ``The Use of Confidence or Fiducial Limits in the Case of the Binomial,'' \textit{Biometrika}, 26, 404-413.
		
DiCiccio, T. J., and Romano, J. P. (1988). ``A Review of Bootstrap Confidence Intervals (with discussions),'' {\it Journal of the Royal Statistical Society B}, 50, 338-354.
	
DiCiccio, T. J., and Tibshirani, R. J. (1987). ``Bootstrap Confidence Intervals
 and Bootstrap Approximations,'' {\it Journal of the American Statistical Association}, 82, 163-170.
		
Efron, B. (1979). ``Bootstrap Methods: Another Look at the Jackknife,'' \textit{Annals of Statistics}, 7, 1-26.
		
Efron, B., Rogosa, D., and Tibshirani, R. (2004). ``Resampling Methods of Estimation,''  In  {\it International Encyclopedia of the Social and Behavioral Sciences} (pp. 13216--13220), Smelser, N. J. and P.B. Baltes, P. B. (Eds.), New York, NY: Elsevier.
		
Efron, B., and Tibshirani, R.J. (1993). \textit{An Introduction to the Bootstrap}. Chapman \& Hall/CRC.
		
Gart, J. J. (1966). ``Alternative Analyses of Contingency Table,''  \textit{Journal of the Royal Statistical Society, Series B}, 28, 164-179.
		
Huwang, L. (1995). ``A Note on the Accuracy of an Approximate Interval for the Binomial Parameter.'' \textit{Statistics and Probability	Letters}, 24, 177-180.
		
Mantalos, P., and Zografos, K. (2008).
	``Interval Estimation for a Binomial Proportion: a Bootstrap Approach.'' \textit{Journal of Statistical Computation and Simulation}, 78, 1251-1265.
	
	Shao, J., and Tu, D. (1995). {\it The Jackknife	and Bootstrap}, Springer-Verlag New York, Inc..
		
Wang, W. (2013). ``A Note on Bootstrap Confidence Intervals for Proportions,''  \textit{Statistics and Probability	Letters}, 83 , 2699-2702.
		
Wang, W. (2014). ``An Iterative Construction of Confidence Interval for a Proportion,''  \textit{Statistica Sinica}, 24, 1389-1410.
		
Wilson, E. B. (1927). ``Probable Inference, the Law of Succession, and Statistical Inference,''  \textit{Journal of the American Statistical Association}, 22, 209-212.
		
Woolf, B. (1955). ``On Estimating the Relation Between Blood Group and Disease,'' \textit{Annals of Human Genetics}, 19, 251-253.

\end{document}